\documentclass[twoside]{article}

\usepackage[margin=2cm]{geometry}
\usepackage{amsmath,amsthm,amssymb,latexsym}
\usepackage[v2,tips]{xy}

\theoremstyle{plain}%
\newtheorem{proposition}{Proposition}[section]%
\newtheorem{theorem}[proposition]{Theorem}%
\newtheorem{corollary}[proposition]{Corollary}%
\newtheorem{lemma}[proposition]{Lemma}%
\theoremstyle{definition}%
\newtheorem{definition}[proposition]{Definition}%

\newcommand{\wap}{\operatorname{WAP}}
\newcommand{\ap}{\operatorname{AP}}

\newcommand{\FIN}{\operatorname{FIN}}

\newcommand{\proten}{{\widehat{\otimes}}}
\newcommand{\mc}[1]{\mathcal{#1}}
\newcommand{\mf}[1]{\mathfrak{#1}}
\newcommand{\ip}[2]{{\langle {#1} , {#2} \rangle}}
\newcommand{\inten}{{\check{\otimes}}}
\newcommand{\aone}{\Box}
\newcommand{\atwo}{\Diamond}

\newcommand{\id}{\operatorname{id}}
\newcommand{\lin}{\operatorname{lin}}
\newcommand{\rank}{\operatorname{rank}}
\newcommand{\Alt}{\operatorname{Alt}}

\begin{document}

\large
\title{Amenability of ultrapowers of Banach algebras}
\author{Matthew Daws}
\maketitle

\begin{abstract}
We study when certain properties of Banach algebras are stable
under ultrapower constructions.  In particular, we consider
when every ultrapower of $\mc A$  is Arens regular, and give some
evidence that this is so if and only if $\mc A$ is isomorphic to
a closed subalgebra of operators on a super-reflexive Banach space.
We show that such ideas are closely related to whether one can
sensibly define an ultrapower of a dual Banach algebra.  We study
how tensor products of ultrapowers behave, and apply this to study
the question of when every ultrapower of $\mc A$ is amenable.
We provide an abstract characterisation in terms of something like
an approximate diagonal, and consider when every ultrapower of a
C$^*$-algebra, or a group $L^1$-convolution algebra, is amenable.


\end{abstract}

\section{Introduction}

Given a Banach space $E$ and an ultrafilter $\mc U$, we can form
the ultrapower $(E)_{\mc U}$.  This construction has proved to
be useful in Banach space theory, especially with regards to local
theory.  Given a Banach algebra $\mc A$, it is trivial that
$(\mc A)_{\mc U}$ is a Banach algebra.  Right at the beginning of
the study of ultrapowers, in \cite{DCK}, certain sequence spaces
which are Banach algebras were studied.  As noted in \cite{Hein}, 
C$^*$-algebra techniques can be used to show that the class of $C(K)$
spaces is closed under ultrapower constructions.  In \cite{GI},
\cite{IL} and \cite{Daws}, ultrapowers of Banach algebras were
used to study the Arens products on the bidual of $\mc A$.  In
\cite{GH}, Ge and Hadwin study ultrapowers of C$^*$-algebras.
Otherwise, the study of ultrapowers
of Banach algebras has been surprisingly sparse
(see below for further points).

For a property of Banach spaces $(P)$, we say that a Banach space
has \emph{super-$(P)$} if every ultrapower of $E$ has $(P)$.
The best known example is that of a \emph{super-reflexive} Banach
space (see \cite[Section~6]{Hein}).  We shall study some super
properties of Banach algebras: in particular, when ultrapowers
of a Banach algebra are Arens regular, and when they are amenable.

There seems to be a close relationship between a Banach algebra
being super Arens regular, and the algebra being isomorphic to
a closed subalgebra of operators on a super-reflexive Banach space.
We also show that the natural construction of an ultrapower of
a dual Banach algebra only works, in practice, for super
Arens regular Banach algebras.

We say that a Banach algebra $\mc A$ is \emph{ultra-amenable}
if every ultrapower of $\mc A$ is amenable (the term
\emph{super-amenable} is used for another meaning by Runde in
\cite{RundeBook}).  We show that ultra-amenability is strictly
weaker than contractability (which is what Runde calls
super-amenable), and strictly stronger than amenability.  Part of
our motivation is that it is generally easy to show that a Banach
algebra is not contractible, while amenability is a much harder
property to settle (this applies in particular to $\mc B(E)$,
the algebra of operators on a Banach space $E$).  We hope that perhaps
the ultra-amenability of $\mc B(E)$ can be more easily settled,
although our current techniques do not allow this.

We provide an abstract characterisation of ultra-amenability,
similar to the concept of an \emph{approximate diagonal},
see \cite[Section~2.2]{RundeBook}.  To do this, we need to first
study how tensor products and ultrapowers interact.  We present
a counter-example due to Charles Read that ultrapowers and
tensor products do not ``commute''.  We settle
when a C$^*$-algebra is ultra-amenable, and show for many
locally compact groups $G$ that ultra-amenability is
equivalent to being finite.

\subsection{Notation and basic concepts}\label{Intro}

We generally follow \cite{Dales} for notation and Banach algebra
concepts.
Let $E$ be a Banach space.  We write $E'$ for the dual space of $E$,
and for $x\in E$ and $\mu\in E'$, we write $\ip{\mu}{x}$ for $\mu(x)$.
We occasionally use square brackets for inner-products.  Recall the
canonical map $\kappa_E:E\rightarrow E''$ defined by
$\ip{\kappa_E(x)}{\mu} = \ip{\mu}{x}$ for $x\in E$ and $\mu\in E'$.
When $\kappa_E$ is an isomorphism, we say that $E$ is \emph{reflexive}.

Recall the notions of filter and ultrafilter.  Let $\mc U$ be a
non-principal ultrafilter on a set $I$, and let $E$ be a Banach
space.  We form the Banach space
\[ \ell^\infty(E,I) = \Big\{ (x_i)_{i\in I}\subseteq E : \|(x_i)\| :=
\sup_{i\in I} \|x_i\| < \infty \Big\}, \]
and define the closed subspace
\[ \mc N_{\mc U} = \Big\{ (x_i)_{i\in I} \in \ell^\infty(E,I) :
\lim_{i\rightarrow\mc U} \|x_i\| =0 \Big\}. \]
Thus we can form the quotient space, called
the \emph{ultrapower of $E$ with respect to $\mc U$},
\[ (E)_{\mc U} := \ell^\infty(E,I) / \mc N_{\mc U}. \]
In general, this space will depend on $\mc U$, though many properties
of $(E)_{\mc U}$ turn out to be independent of $\mc U$, as long as
$\mc U$ is sufficiently ``large'' in some sense.

We can verify that, if $(x_i)_{i \in I}$ represents an equivalence
class in $(E)_{\mc U}$, then
\[ \| (x_i)_{i\in I} + \mc N_{\mc U} \|
= \lim_{i\rightarrow\mc U} \|x_i\|. \]
We shall abuse notation and write $(x_i)$ for the equivalence
class it represents; of course, it can be checked that any
definition we make is independent of the choice of representative
of the equivalence class.  There is a canonical isometry
$E \rightarrow (E)_{\mc U}$ given by sending $x\in E$ to
the constant family $(x)$.  We again abuse notation and write
$x \in (E)_{\mc U}$, identifying $E$ with a closed subspace of $(E)_{\mc U}$.

\begin{definition}
An ultrafilter $\mc U$ is \emph{countably incomplete}
when there exists a sequence $(U_n)_{n=1}^\infty$ in $\mc U$
such that $U_1 \supseteq U_2 \supseteq U_3 \supseteq \cdots$
and such that $\bigcap_n U_n=\emptyset$.
\end{definition}

Countably incomplete ultrafilters are useful, because they allow
us to embed sequential convergence into convergence along the
ultrafilter (see numerous examples of this argument in \cite{Hein}).
We remark that if there exists a non-countably incomplete ultrafilter,
then there exists an uncountable \emph{measurable cardinal}, and it is known that the
existence of such cardinals cannot be shown in ZFC.  See \cite[Section~4.2]{CK}
for further details.  Notice that any non-principal ultrafilter on a countable index
set is certainly countably incomplete.

There is a canonical map $(E')_{\mc U} \rightarrow (E)_{\mc U}'$ given by
\[ \ip{(\mu_i)}{(x_i)} = \lim_{i\rightarrow\mc U} \ip{\mu_i}{x_i}
\qquad ( (\mu_i)\in (E')_{\mc U}, (x_i)\in (E)_{\mc U} ). \]
This map is an isometry, and so we identify $(E')_{\mc U}$ with
a closed subspace of $(E)'_{\mc U}$.  It is shown in \cite[Proposition~7.1]{Hein} that
when $\mc U$ is countably incomplete, $(E)'_{\mc U} = (E')_{\mc U}$ if and only if
$(E)_{\mc U}$ is reflexive.  Furthermore, we define a Banach space $E$ to
be \emph{super-reflexive} if $(E)_{\mc U}$ is reflexive for any ultrafilter
$\mc U$.  As shown in \cite[Proposition~6.4]{Hein}, this definition is equivalent
to the original one given by James (see \cite{James1}).

For Banach spaces $E$ and $F$, we write $\mc B(E,F)$
for the space of bounded linear operators from $E$ to $F$.  Then
there is a canonical isometric map $( \mc B(E,F) )_{\mc U}
\hookrightarrow \mc B( (E)_{\mc U}, (F)_{\mc U} )$ given by
\[ T(x) = (T_i(x_i)) \qquad ( T=(T_i)\in ( \mc B(E,F) )_{\mc U},
x=(x_i)\in (E)_{\mc U} ). \]
We shall often identify $(\mc B(E,F))_{\mc U}$ with its image
in $\mc B( (E)_{\mc U}, (F)_{\mc U} )$.

\section{Basics of ultrapowers of Banach algebras}

When $\mc A$ is a Banach algebra, $(\mc A)_{\mc U}$ becomes
a Banach algebra under the pointwise product.  This follows, as it
is easy to show that $\mc N_{\mc U}$ is a closed ideal in the Banach
algebra $\ell^\infty(\mc A,I)$.

In \cite{GH}, Ge and Hadwin make a general study of ultrapowers of
C$^*$-algebras.  Much of what they prove can easily be adapted to
general Banach algebras.  To give just one example, the ideas of
\cite[Section~3]{GH} will show that if $\mc A$ is a separable
Banach algebra and $\mc U$ and $\mc V$ are non-principal ultrafilters
on $\mathbb N$, then $(\mc A)_{\mc U}$ and $(\mc A)_{\mc V}$ are
isomorphic as Banach algebras, assuming the continuum hypothesis holds.

\begin{proposition}
For a Banach algebra $\mc A$, an ultrapower $(\mc A)_{\mc U}$
is unital if and only if $\mc A$ is unital.
\end{proposition}
\begin{proof}
Clearly, if $\mc A$ is unital, then so is $(\mc A)_{\mc U}$.
Let $e = (e_i) \in (\mc A)_{\mc U}$ be a unit for $(\mc A)_{\mc U}$,
and choose $e_i$ such that $\|e_i\| = \|e\|\geq 1$ for each $i$.
For $\epsilon>0$, for each $i$, let $a_i\in\mc A$ be such that $\|a_i\|=1$ and
\[ \|a_i - e_ia_i\| \geq
\sup\{ \|a-e_ia\| : a\in\mc A, \|a\|\leq1 \} - \epsilon. \]
Let $b=(a_i)\in(\mc A)_{\mc U}$, so that $b = eb = (e_ia_i)
\in(\mc A)_{\mc U}$, and hence
\[ 0 = \lim_{i\rightarrow\mc U} \|a_i - e_ia_i\| \geq
\lim_{i\rightarrow\mc U} 
\sup\{ \|a-e_ia\| : a\in\mc A, \|a\|\leq1 \} - \epsilon. \]
As $\epsilon>0$ was arbitrary, we see that
\[ \lim_{i\rightarrow\mc U}
\sup\{ \|a-e_ia\| : a\in\mc A, \|a\|\leq1 \} = 0. \]
Analogously, we see that
\[ \lim_{i\rightarrow\mc U}
\sup\{ \|a - a e_i\| : a\in\mc A, \|a\|\leq1 \} = 0. \]

For $\epsilon>0$, let $U = \{i : \|a-ae_i\| + \|a-e_ia\|<\epsilon \
(a\in\mc A, \|a\|\leq 1) \} \in\mc U$.  Thus, for $i,j\in U$,
\[ \| e_i - e_j \| \leq \| e_i - e_ie_j \| + \| e_j - e_ie_j \|
< 2\epsilon\|e\|. \]
It is straightforward to extract a sequence from the family
$(e_i)$ which will be Cauchy, and hence converges to,
say, $e_{\mc A}\in\mc A$.  It is then clear that $e_{\mc A}$ will be
a unit for $\mc A$.
\end{proof}

The following is perhaps a little more surprising.

\begin{proposition}\label{bai_ultra}
For a Banach algebra $\mc A$, an ultrapower $(\mc A)_{\mc U}$
has a bounded approximate identity if and only if $\mc A$ does.
The same statement holds for left or right bounded approximate
identities.
\end{proposition}
\begin{proof}
Suppose that $\mc A$ has a bounded approximate identity of bound $M\geq1$.
Let $\mc U$ be an ultrafilter on an index set $I$, let $a=(a_i)\in(\mc A)_{\mc U}$,
and let $\epsilon>0$.  Again, we may suppose that $\|a_i\| = \|a\|$ for
each $i\in I$.  For each $i\in I$, we can find $u_i\in\mc A$ with $\|u_i\|\leq M$
and $\|a_i-u_ia_i\|<\epsilon$.  Let $u=(u_i)_{i\in I}\in (\mc A)_{\mc U}$, so
that $\|a-ua\|\leq\epsilon$.  As $a$ and $\epsilon$ were arbitrary, we see that
$(\mc A)_{\mc U}$ has \emph{bounded left approximate units}.
By \cite[Corollary~2.9.15]{Dales}, we have that $(\mc A)_{\mc U}$
has a bounded left approximate identity of bound $M$.  By symmetry,
$(\mc A)_{\mc U}$ has a bounded right approximate identity of bounded $M$,
and so by a result due to Dixon, see \cite[Proposition~2.9.3]{Dales},
we have that $(\mc A)_{\mc U}$ has a bounded approximate identity
of bound $2M+M^2$.

Conversely, suppose that $(\mc A)_{\mc U}$ has a bounded approximate
identity of bound $M$, but that $\mc A$ does not have a bounded left
approximate identity of bound $\leq M$.  Hence $\mc A$ does not have
bounded left approximate units of bound $\leq M$.  In particular, there
exists $a\in\mc A$ and $\delta>0$ such that $\|a-ua\|\geq\delta$ for
all $u\in\mc A$ with $\|u\|\leq M$.  However, we can find
$u=(u_i)\in(\mc A)_{\mc U}$ with $\| a-ua \| < \delta/2$, that
is, $\lim_{i\rightarrow\mc U} \| a - u_i a \| < \delta/2$, a contradiction.
So $\mc A$ has a bounded left approximate identity of bound $M$, and
thus by symmetry, $\mc A$ has a bounded approximate identity.
\end{proof}

Ultrapowers have been studied in the context of von Neumann algebras.
However, here the definition is different to ours: this is because, for
example, if $\mc M=\ell^\infty$ and $\mc U$ is a non-principal ultrafilter on
$\mathbb N$, then $(\mc M)_{\mc U}$ is not a dual space, and hence not
a von Neumann algebra.  Instead, a construction using traces is often used;
however, it can be shown that the predual of the von Neumann algebra ultrapower is
precisely the Banach space ultrapower of the predual (see, for example,
\cite[Section~1]{raynaud}).  We study such ideas for dual Banach algebras below.

Ultrapowers of Banach spaces have been used in \cite{CF} to study
representations of Banach algebras and representations of groups;
see also the similar ideas used in \cite{DawsArens} and \cite{RundeRep}.

\section{Arens regularity}

Let $\mc A$ be a Banach algebra.  We now recall the Arens products
on $\mc A''$.  Firstly, we turn $\mc A'$ into a $\mc A$-bimodule in
the usual fashion,
\[ \ip{a\cdot\mu}{b} = \ip{\mu}{ba}, \quad
\ip{\mu\cdot a}{b} = \ip{\mu}{ab} \qquad (a,b\in\mc A,\mu\in\mc A'). \]
In a similar way, $\mc A''$ and so forth also become $\mc A$-bimodules.
Then we define bilinear maps $\mc A''\times\mc A', \mc A'\times\mc A''\rightarrow
\mc A'$ by
\[ \ip{\Phi\cdot\mu}{a} = \ip{\Phi}{\mu\cdot a}, \quad
\ip{\mu\cdot\Phi}{a} = \ip{\Phi}{a\cdot\mu} \qquad (\Phi\in\mc A'',
\mu\in\mc A',a\in\mc A). \]
Finally, we define bilinear maps $\aone,\atwo:\mc A''\times\mc A''\rightarrow
\mc A''$ by
\[ \ip{\Phi\aone\Psi}{\mu} = \ip{\Phi}{\Psi\cdot\mu},
\quad \ip{\Phi\atwo\Psi}{\mu} = \ip{\Psi}{\mu\cdot\Phi}
\qquad (\Phi,\Psi\in\mc A'', \mu\in\mc A'). \]
These are associative products which extend the natural action
of $\mc A$ on $\mc A''$, called the \emph{first} and \emph{second Arens products}.
See \cite[Section~3.3]{Dales} or \cite[Section~1.4]{Palmer1} for further details.
Thus $\aone$ and $\atwo$ agree with the usual product on $\kappa_{\mc A}(\mc A)$.
When $\aone$ and $\atwo$ agree on all of $\mc A''$, we say that
$\mc A$ is \emph{Arens regular}.

By Goldstein's Theorem, we know that the unit ball of $\mc A$ is weak$^*$-dense in
the unit ball of $\mc A''$.  This allows us to find an ultrafilter $\mc U$ such
that, given $\Phi,\Psi\in\mc A''$, we can find bounded families $(a_i)$ and $(b_i)$
with $(a_i)$ tending to $\Phi$ weak$^*$ along $\mc U$, and $(b_i)$ tending to $\Psi$.
See \cite[Proposition~6.7]{Hein} for further details.  Then
\[ \ip{\Phi\aone\Psi}{\mu} =
\lim_{j\rightarrow\mc U} \lim_{i\rightarrow\mc U} \ip{\mu}{a_ib_j}, \quad
\ip{\Phi\atwo\Psi}{\mu} =
\lim_{i\rightarrow\mc U} \lim_{j\rightarrow\mc U} \ip{\mu}{a_ib_j}
\qquad (\mu\in\mc A'). \]
In \cite{Daws} we show that when $\mc A$ is Arens regular, we can find
a more ``symmetric'' version of these formulae.

We shall say that $\mc A$ is \emph{super Arens regular} if every
ultrapower of $\mc A$ is Arens regular.  As Arens regularity passes to
subalgebras, clearly a super Arens regular Banach algebra is Arens regular.

\begin{proposition}
Let $\mc A$ be a Banach algebra isomorphic to a closed subalgebra of
$\mc B(E)$ for a super-reflexive Banach space $E$.  Then $\mc A$ is
super Arens regular.
\end{proposition}
\begin{proof}
It is shown in \cite{DawsArens} that $\mc B(E)$ is Arens regular for
any super-reflexive Banach space $E$.  Let $\mc U$ be an ultrafilter.
As an ultrapower of an ultrapower is again an ultrapower (see
\cite[Page~90]{Hein}), we see that $(E)_{\mc U}$ is super-reflexive.
We identify $(\mc B(E))_{\mc U}$ as a closed subalgebra of
$\mc B( (E)_{\mc U} )$.  Thus $(\mc B(E))_{\mc U}$ is Arens regular,
and hence so is $(\mc A)_{\mc U}$, as required.
\end{proof}

Note that if $\mc A$ is a Banach algebra whose underlying
Banach space is super-reflexive, then every ultrapower of $\mc A$
is reflexive, and hence certainly Arens regular.  As noted in \cite{DawsArens},
if $\mc A$ is a closed subalgebra of $\mc B(E)$ for a super-reflexive $E$
(or $\mc A$ is super-reflexive) then every even dual of $\mc A$ is
Arens regular.

Let $\mu\in\mc A'$.  We say that $\mu$ is \emph{weakly almost periodic}
if the map
\[ L_\mu: \mc A\rightarrow\mc A'; \quad a\mapsto a\cdot\mu \qquad (a\in\mc A), \]
is weakly compact, and write $\mu\in\wap(\mc A')$.  Then $\mc A$ is
Arens regular if and only if $\wap(\mc A') = \mc A'$.  See
\cite[Section~3]{DL} for further details (and be aware that they write
$\wap(\mc A)$).  A useful characterisation of $\wap(\mc A')$, due originally
to John Pym (see \cite[Theorem~4.3]{pym}), is that $\mu\in\wap(\mc A')$ if and only
if $\ip{\Phi\aone\Psi}{\mu} = \ip{\Phi\atwo\Psi}{\mu}$ for all
$\Phi,\Psi\in\mc A''$.  Combining this fact with some careful
arguments yields the following repeated limit criterion.

\begin{proposition}
Let $\mc A$ be a Banach algebra, and let $\mu\in\mc A'$.  Then
$\mu$ is weakly almost periodic if and only if, for bounded
sequences $(a_n)$ and $(b_m)$ in $\mc A$, we have that
\[ \lim_{n\rightarrow\infty} \lim_{m\rightarrow\infty} \ip{\mu}{a_n b_m}
= \lim_{m\rightarrow\infty} \lim_{n\rightarrow\infty} \ip{\mu}{a_n b_m}, \]
whenever all the iterated limits exist.
\end{proposition}
\begin{proof}
See \cite[Theorem~2.6.17]{Dales} or \cite[Section~3]{DL}, for example.
\end{proof}

For an ultrapower $(\mc A)_{\mc U}$, we generally do not fully
understand the dual $(\mc A)_{\mc U}'$.  However, we have the norming
subspace $(\mc A')_{\mc U}$, and so in particular, if $(\mc A)_{\mc U}$
is Arens regular, then $(\mc A')_{\mc U} \subseteq \wap((\mc A)_{\mc U}')$.

\begin{lemma}
Let $\mc A$ be a Banach algebra, and let $\mc U$ be a countably incomplete ultrafilter.
An ultrapower $(\mc A)_{\mc U}$ is
Arens regular if and only if $(\mc A')_{\mc U} \subseteq \wap((\mc A)_{\mc U}')$.
\end{lemma}
\begin{proof}
We need only show the ``if'' part.  Let $\mu\in(\mc A)_{\mc U}'$, and
suppose that $\mu$ is not weakly almost periodic.  Thus there exist
bounded sequences $(a_n)$ and $(b_m)$ in $(\mc A)_{\mc U}$ such that
the iterated limits of $(\ip{\mu}{a_nb_m})$ exist, but are not equal.
Let $E\subseteq(\mc A)_{\mc U}$ be the closed linear span of
$(a_nb_m)_{n,m\in\mathbb N}$, so that $E$ is separable.
As $\mc U$ is countably incomplete, we can apply
\cite[Corollary~7.5]{Hein} to see that there exists $\lambda\in(\mc A')_{\mc U}$
such that $\ip{\mu}{a_nb_m} = \ip{\lambda}{a_nb_m}$ for all $n$ and $m$.
Thus $\lambda$ is not weakly almost periodic, a contradiction.
\end{proof}

For a Banach space $E$, let $\mc F(E)$ be the space of finite-rank
operators on $E$, and let $\mc A(E)$ be the space of \emph{approximable
operators}, the norm closure of $\mc F(E)$ in $\mc B(E)$.  See below
for further details, and for what it means for a Banach space $E$ to
have the \emph{approximation property}.  It is known that if $E$ is
a reflexive Banach space with the approximation property, then
$\mc A(E)$ is Arens regular, and that $\mc A(E)'' = \mc B(E)$ as
a Banach algebra (see, for example, \cite{PalmerBiDual} or
\cite[Section~1.7]{Palmer1}).  In general, $\mc A(E)$ is Arens regular
if and only if $E$ is reflexive (see \cite[Theorem~2.6.23]{Dales}, for example).

\begin{proposition}
Let $E$ be a Banach space.  The $\mc A(E)$ is super Arens regular
if and only if $E$ is super-reflexive.
\end{proposition}
\begin{proof}
By the above, if $E$ is super-reflexive, then $\mc A(E)$ is super Arens regular.
If $E$ is not super-reflexive, then by the results of \cite[Section~6]{Hein},
there exists a countably incomplete ultrafilter $\mc U$ on an index set $I$ such
that $(E)_{\mc U}$ is not reflexive.  By a result of James (see, for example,
\cite[Section~4]{DawsArens}), we can find bounded sequences $(x^{(n)})$
in $(E)_{\mc U}$ and $(\mu^{(m)})$ in $(E)_{\mc U}'$ such that
\[ \ip{\mu^{(m)}}{x^{(n)}} = \begin{cases} 0 &: m>n, \\ 1 &:m\leq n.
\end{cases} \]
Let $E$ be the closed linear span of the $(x^{(n)})$, so that $E$
is separable.  As we only care about the value of $\mu^{(m)}$ on $E$,
by \cite[Corollary~7.5]{Hein}, we may suppose that $\mu^{(m)}\in(E')_{\mc U}$.
Let $x^{(n)}=(x^{(n)}_i)$ and $\mu^{(n)} = (\mu^{(n)}_i)$ for each $n$.

Let $\lambda\in E'$ and $x\in E$ be such that $\ip{\lambda}{x}=1$.
For each $n\geq 1$, define
\[ T_n = \mu^{(n)} \otimes x = (\mu^{(n)}_i\otimes x)\in(\mc A(E))_{\mc U}, \quad
S_n = \lambda\otimes x^{(n)} = (\lambda\otimes x^{(n)}_i)\in(\mc A(E))_{\mc U}. \]
Define $\Lambda\in(\mc A(E))_{\mc U}'$ by
\[ \ip{\Lambda}{R} = \lim_{i\rightarrow\mc U} \ip{\lambda}{R_i(x)}
\qquad (R=(R_i)\in (\mc A(E))_{\mc U}). \]
It is hence easy to see that
\[ \ip{\Lambda}{T_n S_m} = \ip{\mu^{(n)}}{x^{(m)}}, \]
from which it follows that $\Lambda$ is not weakly almost periodic,
as required.
\end{proof}

Notice that in the above proof, $\Lambda$ is a member of $\mc A(E)'$,
where we naturally embed $\mc A(E)'$ into $(\mc A(E))_{\mc U}'$.  Hence,
when $E$ is not super-reflexive, $(\mc A(E))_{\mc U}$ fails to be
Arens regular is this rather strong sense.

An alternative way to see the above is the following.  For a Banach
space $E$, we can regards $(\mc A(E))_{\mc U}$ as a subalgebra of
$\mc B( (E)_{\mc U} )$ in the usual way.  It is then easy to see that
$\mc A( (E)_{\mc U} )$ is contained in $(\mc A(E))_{\mc U}$, and so if
$(\mc A(E))_{\mc U}$ is Arens regular, so is $\mc A( (E)_{\mc U} )$,
and hence, as mentioned above, $(E)_{\mc U}$ must be reflexive.

Combining the above results, we might be tempted to make the
following conjecture: a Banach algebra $\mc A$ is super Arens regular
if and only if $\mc A$ is isomorphic to a subalgebra of $\mc B(E)$
for some super-reflexive Banach space $E$.  In \cite{Young},
Young showed that a Banach algebra $\mc A$ is isomorphic to
a subalgebra of $\mc B(E)$ for a reflexive $E$ if and only if
$\wap(\mc A')$ approximately norms $\mc A$, that is, for some
$\delta>0$,
\[ \|a\| \geq \delta \sup\big\{ |\ip{\mu}{a}| : \mu\in\wap(\mc A'),
\|\mu\|\leq 1 \big\} \qquad (a\in\mc A). \]
In particular, Arens regular Banach algebras are even isometric
to closed subalgebras of $\mc B(E)$ for reflexive $E$.
The key tool which Young uses is that of interpolation spaces,
although this wasn't recognised at the time (compare Kaijser's work
in \cite{Kai}).  However, it is not clear how interpolation spaces
and ultrapowers interact; just because an ultrapower $(\mc A)_{\mc U}$
is isomorphic to a subalgebra of $\mc B(E)$ does not seem to
imply that $E$ need to be ultrapower.

\subsection{Ultrapowers of dual Banach algebras}\label{ultra_dual_ba}

Surprisingly, defining ultrapowers of dual Banach algebras is
not as straight forward as for von Neumann algebras: we have to
take account of Arens regularity.

Recall that a \emph{dual Banach algebra} is a Banach algebra
$\mc A$ which is the dual of a Banach space, say $\mc A = \mc A_*'$,
such that the product on $\mc A$ is separately weak$^*$-continuous.
The canonical example is a von Neumann algebra, in which case
the predual $\mc A_*$ is isometrically unique.  In general, there
may be a choice of $\mc A_*$, so we shall write $(\mc A,\mc A_*)$
to indicate the predual.  See \cite{Runde} or
\cite{DawsDBA} for general further information.

By analogy with the von Neumann case, the natural way to define
an \emph{ultrapower} of $\mc A$ is to form the Banach space
ultrapower $(\mc A_*)_{\mc U}$, and then to extend the product
from $(\mc A)_{\mc U}$ to the dual space $(\mc A_*)_{\mc U}'$.

\begin{proposition}\label{up_dba}
Let $(\mc A,\mc A_*)$ be a dual Banach algebra, and let $\mc U$
be an ultrafilter on an index set $I$.  Let $\mf A_* = (\mc A_*)_{\mc U}$
and $\mf A = \mf A_*'$.  The following are equivalent:
\begin{enumerate}
\item\label{updba:one} There is a product on $\mf A$ extending the
   product on $(\mc A)_{\mc U}$ and turning $(\mf A,\mf A_*)$ into
   a dual Banach algebra;
\item\label{updba:two} If we identify $(\mc A_*)_{\mc U}$ with a
   subspace of $(\mc A)_{\mc U}'$, we have that $(\mc A_*)_{\mc U}
   \subseteq \wap( (\mc A)_{\mc U}' )$.
\end{enumerate}
\end{proposition}
\begin{proof}
Notice that as $(\mc A)_{\mc U}$ is weak$^*$-dense in $\mf A$,
any product making $(\mf A,\mf A_*)$ into a dual Banach algebra,
and which extends the product on $(\mc A)_{\mc U}$, must be unique.
If (\ref{updba:one}) holds then it is an easy calculation
(see \cite[Section~2]{DawsDBA}) that $\mf A_* \subseteq \wap(\mf A')$.
Condition (\ref{updba:two}) is immediate from this.

Conversely, notice that $\mf A = (\mc A_*)_{\mc U}$ is an
$(\mc A)_{\mc U}$-bimodule, and so $\mf A'$ is also an
$(\mc A)_{\mc U}$-bimodule.  It is obvious that this bimodule structure
extends the product on $(\mc A)_{\mc U}$.  We can hence extend this
bimodule structure to a bilinear map on $\mf A$, either by extending
on the left, or on the right, by weak$^*$-continuity.  Let us check
that these give the same result.  Let $a,b\in\mf A$, so by
\cite[Section~7]{Hein}, there exist bounded nets $(a_\alpha)$ and
$(b_\alpha)$ in $(\mc A)_{\mc U}$, tending to $a$ and $b$ respectively.
For $\mu\in\mf A_*$, we see that
\begin{align*} \lim_\alpha \ip{a_\alpha\cdot b}{\mu}
&= \lim_\alpha \ip{b}{\mu\cdot a_\alpha}
= \lim_\alpha \lim_\beta \ip{b_\beta}{\mu\cdot a_\alpha}
= \lim_\alpha \lim_\beta \ip{a_\alpha b_\beta}{\mu} \\
&= \lim_\beta \lim_\alpha \ip{a_\alpha b_\beta}{\mu}
= \lim_\beta \ip{a \cdot b_\beta}{\mu}.
\end{align*}
We can swap the order of the limits, as $\mu\in\wap( (\mc A)_{\mc U}' )$.
The construction of this product is very similar to the construction
of the Arens products, and checking that our product on $\mf A$ is
associative is similar to the analogous calculation for the Arens products.

Finally, we show that $(\mf A,\mf A_*)$ is a dual Banach algebra,
for which it suffices to check that $\mf A_*$ is an $\mf A$-submodule
of $\mf A'$.  Let $\mu\in\mf A_*$ and $a\in\mf A$, and suppose that
$a\cdot\mu \not\in\mf A_* \subseteq \mf A'$, so there exists
$\Phi\in\mf A''$ annihilating $\mf A_*$ and with $\ip{\Phi}{a\cdot\mu}=1$.
Let $(b_\alpha)$ be a bounded net in $\mf A$ tending to $\Phi$ weak$^*$
in $\mf A''$.  For each $\alpha$, let $(c_{\alpha,\beta})$ be a bounded
net in $(\mc A)_{\mc U}$ tending to $b_\alpha$ weak$^*$ in $\mf A$.
Let $(c_\gamma)$ be a bounded net in $(\mc A)_{\mc U}$ tending to
$a$ in the weak$^*$-topology on $\mf A$.  Then we see that
\begin{align*} 1 &= \ip{\Phi}{a\cdot\mu}
= \lim_\alpha \ip{a\cdot\mu}{b_\alpha}
= \lim_\alpha \ip{b_\alpha a}{\mu}
= \lim_\alpha \lim_\beta \ip{c_{\alpha,\beta} a}{\mu}
= \lim_\alpha \lim_\beta \ip{a}{\mu\cdot c_{\alpha,\beta}} \\
&= \lim_\alpha \lim_\beta \lim_\gamma \ip{c_\gamma}{\mu\cdot c_{\alpha,\beta}}
= \lim_\alpha \lim_\beta \lim_\gamma \ip{c_{\alpha,\beta} c_\gamma}{\mu}
= \lim_\gamma \lim_\alpha \lim_\beta \ip{c_{\alpha,\beta} c_\gamma}{\mu} \\
&= \lim_\gamma \lim_\alpha \ip{b_\alpha}{c_\gamma \cdot \mu}
= 0,
\end{align*}
a contradiction.  Again, we use that $\mu\in\wap( (\mc A)_{\mc U}' )$
to allow us to swap the order of limits.  Hence $a\cdot\mu\in\mf A_*$,
and similarly $\mu\cdot a\in\mf A_*$, as required.
\end{proof}

Notice that if $\mc A$ is super Arens regular, then certainly
condition (2) above always holds.

\begin{proposition}
Let $(\mc A,\mc A_*)$ be a dual Banach algebra, and suppose that
for all ultrafilters $\mc U$, we have that
$(\mc A_*)_{\mc U} \subseteq \wap( (\mc A)_{\mc U}' )$.
Then every even dual of $\mc A$ is Arens regular.
\end{proposition}
\begin{proof}
Firstly we show that $\mc A$ is Arens regular.  Let $\mu\in\mc A'$
and let $(a_n)$ and $(b_m)$ be bounded sequences in $\mc A$ with
the repeated limits $\lim_n \lim_m \ip{\mu}{a_n b_m}$ and
$\lim_m \lim_n \ip{\mu}{a_n b_m}$ existing.  By \cite[Proposition~6.7]{Hein},
for a suitable ultrafilter $\mc U$, there exists $(\mu_i)\in(\mc A_*)_{\mc U}$ with
\[ \lim_{i\rightarrow\mc U} \ip{a}{\mu_i} = \ip{\mu}{a} \qquad (a\in\mc A). \]
As $(\mc A_*)_{\mc U} \subseteq \wap( (\mc A)_{\mc U}' )$, we have that
$(\mu_i) \in \wap( (\mc A)_{\mc U}' )$, and so
\begin{align*} \lim_n \lim_m \ip{\mu}{a_n b_m} &=
\lim_n \lim_m \lim_{i\rightarrow\mc U} \ip{a_n b_m}{\mu_i}
= \lim_n \lim_m \ip{(a_n) (b_m)}{(\mu_i)} \\
&= \lim_m \lim_n \ip{(a_n) (b_m)}{(\mu_i)}
= \lim_m \lim_n \ip{\mu}{a_n b_m}, \end{align*}
as required.

Let $\Lambda\in\mc A'''$ and let $(\Phi_n)$ and $(\Psi_m)$ be bounded
sequences in $\mc A''$ with the repeated limits
$\lim_n \lim_m \ip{\Lambda}{\Phi_n \Psi_m}$ and
$\lim_m \lim_n \ip{\Lambda}{\Phi_n \Psi_m}$ existing.
For an ultrapower $(\mc A)_{\mc U}$, define a map
$\sigma_{\mc U}: (\mc A)_{\mc U}\rightarrow\mc A''$ by
\[ \ip{\sigma_{\mc U}(a)}{\mu} = \lim_{i\rightarrow\mc U}\ip{\mu}{a_i}
\qquad (a = (a_i)\in(\mc A)_{\mc U}). \]
As $\mc A$ is Arens regular, by the main result of \cite{Daws},
there exists an ultrafilter $\mc U$ on an index set $I$, and a map
$K:\mc A''\rightarrow(\mc A)_{\mc U}$, such that $\sigma_{\mc U}\circ K$
is the identity on $\mc A''$, and
\[ \ip{\sigma_{\mc U}( K(\Phi) K(\Psi) )}{\mu} = \ip{\Phi\aone\Psi}{\mu}
\qquad (\mu\in\mc A', \Phi,\Psi\in\mc A''). \]
There exists an ultrafilter $\mc V$ on an index set $J$ such that
$\sigma_{\mc V}:(\mc A')_{\mc V} \rightarrow \mc A'''$ is surjective.

We define (see the end of Section~7 in \cite{Hein}) the ultrafilter
$\mc U\times\mc V$ on $I\times J$ by, for $A\subseteq I \times J$,
setting $A\in\mc U\times\mc V$ if and only if
\[ \big\{ i\in I : \{ j\in J : (i,j)\in A \} \in \mc V \big\} \in \mc U. \]
Then, for a family $(x_{i,j})_{i\in I, j\in J}$ in a compact Hausdorff
space $X$, we have that
\[ \lim_{j\rightarrow\mc V} \lim_{i\rightarrow\mc U} x_{i,j} =
\lim_{(i,j)\rightarrow\mc U\times\mc V} x_{i,j}. \]

For each $n$ let $K(\Phi_n) = (a^{(n)}_i)\in(\mc A)_{\mc U}$, and
let $K(\Psi_n) = (b^{(n)}_i)\in(\mc A)_{\mc U}$.  Let $(\mu_j)\in
(\mc A')_{\mc V}$ be such that $\sigma_{\mc V}((\mu_j)) = \Lambda$.
We then see that, as $(\mu_j) \in \wap( (\mc A)_{\mc U\times\mc V}' )$,
\begin{align*}
\lim_n \lim_m \ip{\Lambda}{\Phi_n \Psi_m}
&= \lim_n \lim_m \lim_{j\rightarrow\mc V} \ip{\Phi_n\Psi_m}{\mu_j}
= \lim_n \lim_m \lim_{j\rightarrow\mc V} \lim_{i\rightarrow\mc U}
   \ip{\mu_j}{a^{(n)}_i b^{(m)}_i} \\
&= \lim_n \lim_m \lim_{(i,j)\rightarrow\mc U\times\mc V} \ip{\mu_j}{a^{(n)}_i b^{(m)}_i} \\
&= \lim_m \lim_n \lim_{(i,j)\rightarrow\mc U\times\mc V} \ip{\mu_j}{a^{(n)}_i b^{(m)}_i}
= \lim_m \lim_n \ip{\Lambda}{\Phi_n \Psi_m}.
\end{align*}
Hence $\mc A''$ is Arens regular.

Repeating this argument allows us to show that every even dual of
$\mc A$ is Arens regular, as claimed.
\end{proof}

Again, it would be interesting to know if, when every even dual of a Banach
algebra $\mc A$ is Arens regular, we have $\mc A$ is super Arens regular?
In conclusion, we see that our approach to ultrapowers of dual Banach algebras
requires a rather strong condition on the underlying algebra, indeed, in
practice, we need $\mc A$ to be a subalgebra of $\mc B(E)$ for
a super-reflexive Banach space $E$.

\section{Tensor products of ultrapowers}\label{tensor_sec}

We shall now sketch the basics of the theory of tensor products of Banach
spaces.  We refer the reader to the books \cite{Ryan} or \cite{DU} for
introductory treatments of this material, or to the book \cite{DF} for
further information.

For Banach spaces $E$ and $F$, let $E \otimes F$ be the algebraic tensor
product of $E$ and $F$.  We define the \emph{projective tensor norm} by
\[ \|\tau\|_\pi = \inf\Big\{ \sum_{i=1}^n \|x_i\| \|y_i\| :
\tau = \sum_{i=1}^n x_i \otimes y_i \Big\}
\qquad ( \tau\in E\otimes F ). \]
The completion of $E\otimes F$ with respect to $\|\cdot\|_\pi$ is
$E\proten F$, the \emph{projective tensor product of $E$ and $F$}.
$E\proten F$ has the universal property that if $T:E\times F\rightarrow G$
is a bounded bilinear map to a Banach space $G$, then there is a unique
bounded linear map $\widehat T:E\proten F\rightarrow G$ such that $\widehat T(x
\otimes y) = T(x,y)$ for $x\in E$ and $y\in F$.  Every member $\tau \in E
\proten F$ can be written as an absolutely convergent sum
$\tau = \sum_{i=1}^\infty x_i \otimes y_i$,
for some sequences $(x_i) \subseteq E$ and $(y_i) \subseteq F$.

Let $\mc F(E,F)$ be the space of finite-rank operators from $E$
to $F$, and let $\mc A(E,F)$ be the space of \emph{approximable} operators
from $E$ to $F$, the norm closure of $\mc F(E,F)$ in $\mc B(E,F)$.
We can embed $E\otimes F$ into $\mc F(E',F)$ by
\[ \Big( \sum_{i=1}^n x_i\otimes y_i \Big): \mu \mapsto
\sum_{i=1}^n \ip{\mu}{x_i} y_i \qquad (\mu\in E'). \]
This induces the \emph{injective tensor norm} $\|\cdot\|_\epsilon$
on $E\otimes F$, whose completion is $E \inten F$.  In particular, we
can identify $\mc A(E,F)$ with $E' \inten F$.

We shall say that norm $\|\cdot\|$ on $E\otimes F$ is a \emph{reasonable crossnorm}
when:
\begin{enumerate}
\item $\| x\otimes y \| = \|x\| \|y\|$ for $x\in E$ and $y\in F$;
\item for $\mu\in E'$ and $\lambda\in F'$, define $\mu\otimes\lambda:
   E\otimes F \rightarrow \mathbb C$ by $\ip{\mu\otimes\lambda}{x\otimes y}
   = \ip{\mu}{x} \ip{\lambda}{y}$ and linearity.  Then the norm of
   $\mu\otimes\lambda$, with respect to $\|\cdot\|$, is $\|\mu\| \|\lambda\|$.
\end{enumerate}
Suppose that for each pair of Banach spaces $(E,F)$, we have an assignment of
a reasonable crossnorm $\|\cdot\|$ on $E\otimes F$.  Then this assignment is a
\emph{uniform crossnorm} when given pairs $(E_1,F_1)$ and $(E_2,F_2)$ of
Banach spaces, for
$T\in\mc B(E_1,E_2), S\in\mc B(F_1,F_2)$, we have that
$\|T\otimes S\| \leq \|T\| \|S\|$ where we treat $T\otimes S$ as a
linear map $E_1\proten F_1 \rightarrow E_2 \proten F_2$ given by
\[ (T\otimes S)(x\otimes y) = T(x) \otimes S(y) \qquad (x\otimes y\in
E_1\otimes F_1), \]
and linearity.
Then $\|\cdot\|_\pi$ and $\|\cdot\|_\epsilon$ are uniform crossnorms.

The projective tensor product is \emph{projective} in the sense that
if $T$ and $S$ are quotient maps (also called metric surjections) then so
is $T\otimes S : E_1\proten F_1 \rightarrow E_2 \proten F_2$.
Similarly, the injective tensor product is \emph{injective} in that,
when $T$ and $S$ are isometries, then so is
$T\otimes S : E_1\inten F_1 \rightarrow E_2 \inten F_2$.  In general,
the projective tensor norm is not injective, and the injective tensor
norm is not projective.  A useful exception to this is that the map
$\kappa_E \otimes \id : E\proten F\rightarrow E''\proten F$ is always an
isometry onto its range.

We identify the dual of $E\proten F$ with $\mc B(E,F')$ by
\[ \ip{T}{x\otimes y} = \ip{T(x)}{y} \qquad \big(
T\in\mc B(E,F'), x\otimes y\in E\proten F \big), \]
and linearity and continuity.  In particular, $E' \inten F'
= \mc A(E,F')$ isometrically embeds into $(E\proten F)'$.
When one of $E$ or $F$ is \emph{finite-dimensional}, we have equality,
$(E\proten F)' = E' \inten F'$.

As the map $E\proten F\rightarrow E\inten F$ is norm-decreasing with
dense range, we see that the adjoint $(E\inten F)' \rightarrow
(E\proten F)' = \mc B(E,F')$ is norm-decreasing and injective.
We hence identify $(E\inten F)'$ with a space of operators $E\rightarrow F'$,
the \emph{integral operators} $\mc I(E,F')$, and we give $\mc I(E,F')$
the dual norm $\|\cdot\|_{\mc I}$, so that $\mc I(E,F') = (E\inten F)'$.
We have a norm-decreasing map $E'\proten F' \rightarrow \mc I(E,F')$.
It is quite a subtle issue as to when this map is bounded below,
an isometry, or when it is surjective.  See \cite{Ryan} or \cite[Section~16]{DF}
for further details.  However, if one of $E$ or $F$
is finite-dimensional, then $(E\inten F)' = \mc I(E,F') = E'\proten F'$.

We say that a Banach space $E$ has the \emph{approximation property}
when the canonical map $E'\proten E\rightarrow E'\inten E = \mc A(E)$
is injective.  See \cite[Chapter~4]{Ryan} or \cite[Chapter~VIII]{DU}
for further details.  For Banach spaces $E$ and $F$ with the approximation
property, we can hence identify $E'\proten F$ as a space of operators
from $E$ to $F$, called the \emph{nuclear operators}, $\mc N(E,F)$.
In general, $\mc N(E,F)$ is merely a quotient of $E'\proten F$, and we
always give $\mc N(E,F)$ the quotient norm.

\subsection{Ultrapowers}

Let $M$ be a finite-dimensional Banach space and let $\mc U$ be an ultrafilter.
By taking a basis, it is easy to see that $(M)_{\mc U} = M$.  It is shown in
\cite[Lemma~7.4]{Hein}, that
\[ (M \inten E)_{\mc U} = M \inten (E)_{\mc U}
\quad,\quad
(M \proten E)_{\mc U} = M \proten (E)_{\mc U} \]
for every Banach space $E$, and every finite-dimensional $M$,
with equality of norms.

For infinite-dimensional Banach spaces, these equalities are no
longer necessarily true.  However, we can make some useful
statements.

Let $E$ and $F$ be Banach spaces.
There is a canonical map $\psi_0 : (E)_{\mc U} \proten (F)_{\mc U}
\rightarrow (E\proten F)_{\mc U}$, defined using the
tensorial property of $\proten$.  Firstly we define $\psi_0:
(E)_{\mc U} \times (F)_{\mc U} \rightarrow (E\proten F)_{\mc U}$ by
\[ \psi_0 ( x, y ) = (x_i\otimes y_i)
\qquad (x=(x_i)\in (E)_{\mc U} , y=(y_i)\in (F)_{\mc U}). \]
Then we have
\[ \| (x_i\otimes y_i) \| = \lim_{i\rightarrow\mc U} \|x_i\otimes y_i\|_\pi
= \lim_{i\rightarrow\mc U} \|x_i\|\|y_i\|
= \Big(\lim_{i\rightarrow\mc U} \|x_i\|\Big)
\Big(\lim_{i\rightarrow\mc U} \|y_i\|\Big) = \|x\|\|y\|, \]
so that $\psi_0$ is well-defined, and is a norm-decreasing
bilinear map.  Thus $\psi_0$ extends to a norm-decreasing map
$\psi_0 : (E)_{\mc U} \proten (F)_{\mc U} \rightarrow (E\proten F)_{\mc U}$.
For $\tau \in (E)_{\mc U} \otimes (F)_{\mc U}$, choose a representative
$\tau = \sum_{k=1}^n x_k \otimes y_k$.  Let, for each $k$,
$x_k = (x^{(k)}_i)\in (E)_{\mc U}$ and $y_k = (y^{(k)}_i)\in (E)_{\mc U}$.
Then we see that
\[ \psi_0(\tau) = \Big( \sum_{k=1}^n x^{(k)}_i \otimes y^{(k)}_i
\Big)_{i\in I} \in ( E\proten F)_{\mc U}. \]

\begin{proposition}\label{image_psi_zero}
Let $E$ and $F$ be Banach spaces, let $\mc U$ be an ultrafilter
on an index set $I$, and let $\tau\in (E\proten F)_{\mc U}$.
Then the following are equivalent:
\begin{enumerate}
\item for some sequence $(\alpha_n)$ of positive reals with
$\sum_n \alpha_n<\infty$, $\tau = (\tau_i)$ admits a
representation of the form
\[ \tau_i = \sum_{k=1}^\infty x^{(i)}_k \otimes y^{(i)}_k
\in E\proten F \qquad (i\in I), \]
where, for each $i$ and $k$, we have that $\|x^{(i)}_k\|
\|y^{(i)}_k\| \leq \alpha_k$;
\item $\tau$ lies in the image of $\psi_0$.
\end{enumerate}
\end{proposition}
\begin{proof}
Suppose that (1) holds.
By rescaling, we may suppose that $\|x^{(i)}_k\| = \|y^{(i)}_k\|
\leq \alpha_k^{1/2}$ for each $i\in I$ and $k\geq1$.  For each
$k\geq1$, let
\[ x_k = (x^{(i)}_k) \in (E)_{\mc U}, \quad
y_k = (y^{(i)}_k) \in (F)_{\mc U}, \]
so that $\|x_k\| \leq \alpha_k^{1/2}$ and
$\|y_k\| \leq \alpha_k^{1/2}$.  We can hence let
\[ \sigma = \sum_{k=1}^\infty x_k \otimes y_k \in
(E)_{\mc U} \proten (F)_{\mc U}, \]
with $\pi(\sigma) \leq \sum_k \alpha_k$.  Let $\sigma_n =
\sum_{k=1}^n x_k \otimes y_k$ so that $\sigma_n\rightarrow\sigma$
in $(E)_{\mc U} \proten (F)_{\mc U}$.  Then
\begin{align*}
\lim_{n\rightarrow\infty} \big\| \psi_0(\sigma_n) - \tau \big\|
&= \lim_{n\rightarrow\infty} \lim_{i\rightarrow\mc U} \Big\|
   \sum_{k=1}^n x^{(i)}_k \otimes y^{(i)}_k - \tau_i \Big\|_\pi \\
&\leq \lim_{n\rightarrow\infty} \lim_{i\rightarrow\mc U}
   \sum_{k=n+1}^\infty \|x^{(i)}_k\| \|y^{(i)}_k\|
\leq \lim_{n\rightarrow\infty} \sum_{k=n+1}^\infty \alpha_k
= 0,
\end{align*}
so that $\psi_0(\sigma) = \tau$, as required.

Conversely, suppose that $\tau = \psi_0(\sigma)$ for
\[ \sigma = \sum_{k=1}^\infty x_k \otimes y_k \in
(E)_{\mc U} \proten (F)_{\mc U}, \]
with $\sum_{k=1}^\infty \|x_k\| \|y_k\|<\infty$.  Then
we let $\alpha_k = \|x_k\| \|y_k\|$ and pick representatives
$x_k = (x^{(i)}_k) \in (E)_{\mc U}$ and
$y_k = (y^{(i)}_k) \in (E)_{\mc U}$, with $\|x_k\| = \|x^{(i)}_k\|$
and $\|y_k\| = \|y^{(i)}_k\|$ for each $k$ and $i$.
For each $i\in I$, let $\tau_i = \sum_{k=1}^\infty x^{(i)}_k
\otimes y^{(i)}_k$.  Let $\sigma_n = \sum_{k=1}^n x_k\otimes
y_k$, so that $\tau = \lim_{n\rightarrow\infty} \psi_0(\sigma_n)$.
Thus, for each $n$,
\begin{align*}
\big\| (\tau_i) - \psi_0(\sigma_n) \big\|
= \lim_{i\rightarrow\mc U} \Big\| \tau_i - \sum_{k=1}^n x^{(i)}_k\otimes
   y^{(i)}_k \Big\|_\pi
\leq \lim_{i\rightarrow\mc U} \sum_{k=n+1}^\infty \|x^{(i)}_k\| \|y^{(i)}_k\|
= \sum_{k=n+1}^\infty \alpha_k.
\end{align*}
Hence, letting $n\rightarrow\infty$, $(\tau_i) = \tau$ as required.
\end{proof}

Let $\mc A$ be a Banach algebra, and let $E$ be a left-$\mc A$-module.
Then an ultrapower $(E)_{\mc U}$ becomes a left-$\mc A$-module in
the obvious way.  When $F$ is a right-$\mc A$-module, we have that
$E\proten F$ is an $\mc A$-bimodule for the module actions
\[ a\cdot(x\otimes y) = a\cdot x \otimes y, \quad
(x\otimes y)\cdot a = x\otimes y\cdot a \qquad (a\in\mc A,
x\otimes y\in E\proten F). \]
Hence an ultrapower $(E\proten F)_{\mc U}$ is also an $\mc A$-bimodule.
Similarly, $(E)_{\mc U} \proten (F)_{\mc U}$ is an $\mc A$-bimodule.
It is a simple check to see that $\psi_0$ is an $\mc A$-bimodule homomorphism.

Similarly, it is easily checked that $(E)_{\mc U}$ is a 
left-$(\mc A)_{\mc U}$-module, $(F)_{\mc U}$ is a right-$(\mc A)_{\mc U}$-module,
and both $(E\proten F)_{\mc U}$ and $(E)_{\mc U} \proten (F)_{\mc U}$
are $(\mc A)_{\mc U}$-bimodules.  We can check that $\psi_0$
is also an $(\mc A)_{\mc U}$-bimodule homomorphism.

In general, it seems that $\psi_0$ is rarely, if ever, surjective
when $E$ and $F$ are infinite-dimensional.  We now present an argument
for Hilbert spaces that is motivated by a counter-example communicated
to us by Charles Read.
We first recall the Schmidt representation
theorem (see, for example, the treatment given in \cite{Pie}).

\begin{theorem}
Let $H$ and $K$ be Hilbert spaces, and $T\in\mc A(H,K)$.
Then there exist orthonormal sequences $(h_n)$ and $(k_n)$
in $H$ and $K$, respectively, and a sequence of positive
numbers $(s_n)\in c_0(\mathbb N)$ with $s_1\geq s_2\geq\cdots$ such that
\[ T(x) = \sum_{n=1}^\infty s_n [x,h_n] k_n
\qquad (x\in H), \]
where $[ \cdot , \cdot ]$ is the inner-product on $H$.
\end{theorem}

Here and henceforth, we allow orthonormal sequences to be eventually zero.

For a Hilbert space $H$ and $x\in H$, we define a linear functional
$x^*$ on $H$ by $y\mapsto [y,x]$.  The Riesz Theorem shows that
every linear functional arises in this way.  It is clear that if
the sequence $(s_n)$ above satisfies $\sum_n s_n < \infty$,
then $T$ will be nuclear, and hence identified with a member of $H\proten K$
(as $H$ and $K$ have the approximation property) with $\|T\|_\pi \leq \sum_n s_n$.

\begin{lemma}
If $T\in H\proten K$, then the sequence $(s_n)$ arising from
the Schmidt representation of $T$ satisfies $\|T\|_\pi = \sum_n s_n$.
\end{lemma}
\begin{proof}
For $\epsilon>0$, let $T = \sum_n u_n^* \otimes v_n$ with
$\sum_n \|u_n\| \|v_n\| < \|T\|_\pi+\epsilon$.  Then, by the
Schmidt representation, we have also that $T = \sum_n s_n
h_n^* \otimes k_n$, say.  Then
\begin{align*}
\sum_n s_n &= \sum_n [ T(h_n), k_n ]
= \sum_n \sum_m [ h_n, u_m] [v_m, k_n] \\
&\leq \sum_m \Big( \sum_n |[ h_n, u_m]|^2 \Big)^{1/2}
\Big( \sum_n |[v_m, k_n]|^2 \Big)^{1/2} \\
&\leq \sum_m \|u_m\| \|v_m\| < \|T\|_\pi + \epsilon,
\end{align*}
as $(h_n)$ and $(k_n)$ are orthonormal sequences.
As $\epsilon>0$ was arbitrary, we are done.
\end{proof}

Notice that this proof shows that, for $T\in H\proten K$,
we have that
\[ \| T \|_\pi = \sup\Big\{ \sum_n |[T(e_n),f_n]| : (e_n)
\text{ and } (f_n) \text{ are orthonormal sequences in }
H \text{ and } K \Big\}. \]

We now recall the notion (see \cite[Chapter~11]{Pie}) of
approximation numbers. Let $E$ and $F$ be Banach spaces, and let
$T\in\mc B(E,F)$.  The \emph{$n$th approximation number of $T$},
for $n\geq 1$, is
\[ a_n(T) = \inf\{ \|T-S\| : S\in\mc F(E,F), \rank(S)<n \}. \]

\begin{proposition}\label{approx_num}
Let $H$ and $K$ be Hilbert spaces, and let $T\in\mc A(H,K)$
have a Schmidt representation $T = \sum_n s_n h_n^*\otimes k_n$.
Then, if $s_1 \geq s_2\geq \cdots$, then $s_n = a_n(T)$ for each
$n\geq 1$.
\end{proposition}
\begin{proof}
See \cite[Section~11.3]{Pie}.
\end{proof}

\begin{proposition}
Let $H$ and $K$ be Hilbert spaces, and let $\mc U$ be a ultrafilter
on an index set $I$.
Then each $\tau \in (H\proten K)_{\mc U}$ in the image of
$\psi_0$ admits a representation of the form $\tau = (\tau_i)$ with
\[ \tau_i = \sum_n s_n h_{n,i}^* \otimes k_{n,i} \in H\proten K
\qquad (i\in I), \]
where $(s_n)$ is a sequence of positive reals with $\sum_n s_n < \infty$ and,
for each $i$, $(h_{n,i})$ and $(k_{n,i})$ are orthonormal sequences in
$H$ and $K$ respectively.
\end{proposition}
\begin{proof}
Let $\tau = \psi_0(\sigma)$ where
\[ \sigma = \sum_n s_n e_n^* \otimes f_n, \]
where $(e_n)$ and $(f_n)$ are orthonormal sequences in
$(H)_{\mc U}$ and $(K)_{\mc U}$ respectively.  Pick representatives
$e_n = (e_{n,i})_{i\in I}$ and $f_n = (f_{n,i})_{i\in I}$, so that
\[ \delta_{n,m} = \lim_{i\rightarrow\mc U} [ e_{n,i} , e_{m,i} ]
= \lim_{i\rightarrow\mc U} [ f_{n,i} , f_{m,i} ] \qquad (n,m\geq 1). \]
For each $i$, apply the Gram-Schmidt orthonormalisation process
to $(e_{n,i})$ to yield $(h_{n,i})$, where we allow $h_{n,i}$ to be
zero for sufficiently large $n$; do the same to $(f_{n,i})$ to
yield $(k_{n,i})$.  For each $n\geq 1$, as $h_{n,i}$ depends
only upon $\{ e_{m,i} : m\leq n \}$ and $\{ [ e_{m,i} , e_{r,i} ]
: m,r\leq n \}$, we can verify that
\[ \lim_{i\rightarrow\mc U} \| h_{n,i} - e_{n,i} \| = 0
\qquad (n\geq 1), \]
and similarly for $f_{n,i}$.  Let $h_n = (h_{n,i})_{i\in I}
\in (H)_{\mc U}$ and $k_n = (k_{n,i})_{i\in I} \in (H)_{\mc U}$,
so that $h_n = e_n$ and $k_n = f_n$.  Thus
\[ \tau = \psi_0(\sigma) = \Big( \sum_n s_n h_{n,i}^* \otimes
k_{n,i} \Big)_{i\in I} \in (H\proten K)_{\mc U}, \]
as required.
\end{proof}

\begin{theorem}
Let $H$ and $K$ be Hilbert spaces, and let $\mc U$ be a countably
incomplete ultrafilter on an index set $I$.  Then $\psi_0:
(H)_{\mc U} \proten (K)_{\mc U} \rightarrow (H\proten K)_{\mc U}$
does not have dense range.
\end{theorem}
\begin{proof}
We first consider the case when $I=\mathbb N$ and $\mc U$ is a
non-principal ultrafilter on $\mathbb N$.  Let $(e_n)$ and $(f_n)$
be orthonormal sequences in $H$ and $K$, respectively.  For each
$n\geq 1$, let
\[ \tau_n = n^{-1} \sum_{j=1}^n e_j^* \otimes f_j \in H\proten K, \]
and let $\tau = (\tau_n) \in (H\proten K)_{\mc U}$.
Let $\sigma\in (H\proten K)_{\mc U}$ be in the image of $\psi_0$,
so that $\sigma$ has a representation as above,
\[ \sigma = (\sigma_k) = \Big( \sum_n s_n h_{n,k}^* \otimes
k_{n,k} \Big)_{k\in\mathbb N} \in (H\proten K)_{\mc U}. \]

Pick $\epsilon>0$, and choose $N$ such that $\sum_{n>N} s_n
< \epsilon$.  Then
\begin{align*}
\lim_{k\rightarrow\mc U} \pi(\tau_k - \sigma_k)
&> \lim_{k\rightarrow\mc U} \pi\Big(\tau_k -
   \sum_{n=1}^N s_n h_{n,k}^*\otimes k_{n,k} \Big) - \epsilon \\
&\geq \lim_{k\rightarrow\mc U} \inf\Big\{ \pi(\tau_k - \upsilon) :
   \upsilon\in\mc F(H,K), \rank(\upsilon)\leq N \Big\} - \epsilon \\
&= \lim_{k\rightarrow\mc U} \inf\Big\{ \sum_{m\geq 1} a_m(\tau_k-\upsilon)
   : \upsilon\in\mc F(H,K), \rank(\upsilon)\leq N \Big\} - \epsilon,
\end{align*}
by an application of Proposition~\ref{approx_num}.  Now,
for $\upsilon\in\mc F(H,K)$ with $\rank(\upsilon)\leq N$, it is
clear that $a_m(\tau_k-\upsilon) \geq a_{N+m}(\tau_k)$, so that
\[ \lim_{k\rightarrow\mc U} \pi(\tau_k - \sigma_k)
> \lim_{k\rightarrow\mc U} (k-N)k^{-1} - \epsilon = 1 - \epsilon. \]
As $\epsilon>0$ and $\sigma$ were arbitrary, we see that $\tau$ is
distance $1$ (as $\pi(\tau)=1$) from the image of $\psi_0$.

A standard argument allows us to adapt this proof to the case when
$\mc U$ is an arbitrary countably incomplete ultrafilter on an
index set $I$ (compare with the proofs of Theorem~6.3 or Proposition~7.1 in \cite{Hein}).
\end{proof}

The above seems to rely very heavily upon certain special features
of Hilbert spaces, as did the original counter-example due to C.\,J.~Read.
It would be interesting to know how $(E)_{\mc U} \proten (F)_{\mc U}$
and $(E\proten F)_{\mc U}$ relate for other classes of Banach spaces.

For the following, we refer the reader to \cite[Section~9]{Hein}, where Heinrich gives a
description of when $(E)_{\mc U}$ has the approximation property.  In particular, the
following are equivalent: $(E)_{\mc U}$ has the approximation property for all $\mc U$;
$E$ has the \emph{uniform approximation property}; $(E)_{\mc U}$ has the approximation
property for some non-principal $\mc U$ on a countable index set.
Notice that, by \cite[Theorem~3.3]{Hein}, $(L^p(\nu))_{\mc U}$ has the approximation
property for any measure $\nu$, $1\leq p\leq\infty$, and any $\mc U$.

\begin{proposition}
Let $E$ and $F$ be Banach spaces such that $F$ is super-reflexive.
Let $\mc U$ be an ultrafilter such that $(F)_{\mc U}$ has the
approximation property.
Then $\psi_0 : (E)_{\mc U} \proten (F)_{\mc U} \rightarrow
(E \proten F)_{\mc U}$ is an isometry onto its range.
\end{proposition}
\begin{proof}
As $(F)_{\mc U}$ is reflexive and $(F)_{\mc U}$ has
the approximation property, we have
\[ \mc A( (E)_{\mc U} , (F')_{\mc U} )'=
( (E)'_{\mc U} \inten (F')_{\mc U} )' =
\mc I( (E)'_{\mc U}, (F)_{\mc U} ) =
(E)''_{\mc U} \proten (F)_{\mc U}. \]
See \cite[Section~5.3]{Ryan} for details.
As the map $\kappa_{(E)_{\mc U}} \otimes \id :
(E)_{\mc U} \proten (F)_{\mc U} \rightarrow
(E)''_{\mc U} \proten (F)_{\mc U}$ is an isometry onto its
range, we see that
\[ \| \tau \|_\pi = \sup\{ |\ip{\tau}{S}| :
S\in\mc F( (E)_{\mc U} , (F')_{\mc U} ), \|S\|\leq1 \}
\qquad (\tau \in (E)_{\mc U} \proten (F)_{\mc U}). \]

In the following, for a Banach space $X$, we write $\FIN(X)$ for
the collection of finite-dimensional subspaces of $X$.
Fix $\tau \in (E)_{\mc U} \otimes (F)_{\mc U}$.
Let $\tau = \sum_{k=1}^n y^{(k)} \otimes z^{(k)}$ and let
$N = \lin\{ y^{(k)} : 1\leq k\leq n\} \in \FIN( (E)_{\mc U} )$.
For each $k$, let $y^{(k)}=(y^{(k)}_i)$ and $z^{(k)}=(z^{(k)}_i)$
where $\|y^{(k)}\|=\|y^{(k)}_i\|$ and $\|z^{(k)}\|=\|z^{(k)}_i\|$
for each $i$.  Thus
\[ \psi_0(\tau) = (\tau_i) = \Big( \sum_{k=1}^n y^{(k)}_i
\otimes z^{(k)}_i \Big)_{i\in I} \in (E\proten F)_{\mc U}. \]

Choose $\epsilon>0$ and let $S\in\mc F( (E)_{\mc U} , (F')_{\mc U} )$
be such that $\|S\|\leq 1$ and $|\ip{\tau}{S}| > \|\tau\|_\pi - \epsilon$.
Let $M = S( (E)_{\mc U} ) \in \FIN( (F')_{\mc U} )$ have a basis
$\{ x^{(1)}, \ldots, x^{(m)} \}$ where $x^{(k)}=(x^{(k)}_i) \in (E)_{\mc U}' =
(E')_{\mc U}$ for each $k$.  Following the proof of \cite[Proposition~6.2]{Hein},
let $M_i=\lin\{ x^{(k)}_i : 1\leq k\leq m\} \in \FIN(F')$ and
$T_i:M \rightarrow M_i$ be defined by $T_i(x^{(k)})=x^{(k)}_i$.  Then, for
some $I_0\in\mc U$, $T_i$ is a $(1+\epsilon)$-isomorphism for each $i\in I_0$.

We can write $S=\sum_{k=1}^m \mu^{(k)}\otimes x^{(k)}$ for some
$(\mu^{(k)})_{k=1}^m \subseteq (E)'_{\mc U}$.  Let $P=\lin\{ \mu^{(k)} : 1\leq k\leq m\}
\in \FIN( (E)_{\mc U}' )$.  By \cite[Theorem~7.3]{Hein}, there exists a
$(1+\epsilon)$-isomorphism onto its range $T : P \rightarrow (E')_{\mc U}$ such that
\[ \ip{T(\mu^{(k)})}{z} = \ip{\mu^{(k)}}{z} \qquad (1\leq k\leq m, z\in N). \]
For each $k$, let $T(\mu^{(k)}) = (\mu^{(k)}_i)\in (E')_{\mc U}$.
Then let $Q=T(P)$, let $Q_i=\lin\{ \mu^{(k)}_i : 1\leq k\leq m \}\in\FIN(E')$
and let $R_i:Q\rightarrow Q_i$ be given by $R_i(T(\mu^{(k)}))=\mu^{(k)}_i$.
Again, there exists $I_1\in\mc U$ such that $R_i$ is
a $(1+\epsilon)$-isomorphism for each $i\in I_1$.

For each $i\in I_0 \cap I_1 \in \mc U$, let
\[ S_i = \sum_{k=1}^m R_iT(\mu^{(k)}) \otimes T_i(x^{(k)})
= (R_iT\otimes T_i)S
\in F' \otimes E' = \mc F(F,E'), \]
so that $\|S_i\| \leq \|R_i\|\|T\|\|T_i\|\|S\|
\leq (1+\epsilon)^3$.
Then we have
\begin{align*}
\ip{\tau}{S} &= \sum_{j=1}^n \ip{S(y^{(j)})}{z^{(j)}}
= \sum_{j=1}^n \sum_{k=1}^m \ip{\mu^{(k)}}{y^{(j)}}
   \ip{x^{(k)}}{z^{(j)}} \\
&= \sum_{j=1}^n \sum_{k=1}^m \ip{T(\mu^{(k)})}{y^{(j)}}
   \ip{x^{(k)}}{z^{(j)}}
= \sum_{j=1}^n \sum_{k=1}^m
   \lim_{i\rightarrow\mc U} \ip{\mu^{(k)}_i}{y^{(j)}_i}
   \ip{x^{(k)}_i}{z^{(j)}_i} \\
&= \lim_{i\rightarrow\mc U} \sum_{j=1}^n \sum_{k=1}^m
   \ip{ R_iT(\mu^{(k)}) }{ y^{(j)}_i }
   \ip{ T_i(x^{(k)}) }{ z^{(j)}_i }
= \lim_{i\rightarrow\mc U} \ip{ S_i }{ \tau_i }.
\end{align*}
As $\psi_0$ is norm-decreasing, we conclude that
\[ \|\tau\|_\pi - \epsilon < \lim_{i\rightarrow\mc U} |\ip{S_i}{\tau_i}|
\leq \lim_{i\rightarrow\mc U} \|S_i\| \pi(\tau_i)
\leq (1+\epsilon)^3 \|\psi_0(\tau)\|
\leq (1+\epsilon)^3 \| \tau \|_\pi. \]
As $\epsilon>0$ was arbitrary, we conclude that $\psi_0$ is
an isometry onto its range.
\end{proof}

\section{Ultra-amenability}

Let $\mc A$ be a Banach algebra and $E$ be a Banach
$\mc A$-bimodule.  Then a \emph{derivation} $d:\mc A
\rightarrow E$ is a linear map such that $d(ab)=a \cdot d(b)
+ d(a)\cdot b$ for $a,b\in\mc A$.  All the derivations which
we shall consider will be continuous.  Let $x\in E$ and
define $\delta_x : \mc A\rightarrow E$ by $\delta_x(a)=
a\cdot x - x\cdot a$.  Then $\delta_x$ is a derivation,
termed an \emph{inner derivation}.

A Banach algebra $\mc A$ is \emph{contractible} or
\emph{super-amenable} if every derivation from $\mc A$
to a Banach $\mc A$-bimodule is inner.  A contractible
Banach algebra is unital, and it is conjectured that
a Banach algebra $\mc A$ is contractible only when $\mc A$
is finite-dimensional.  This is true for C$^*$-algebras
(indeed, for closed subalgebras of $\mc B(H)$ for a Hilbert
space $H$, see \cite{paulsmith}) and for $\mc B(E)$ when $E$ has,
for example, the approximation property (see \cite[Section~4.1]{RundeBook}).

A Banach algebra $\mc A$ is \emph{amenable} if every
derivation from $\mc A$ to a \emph{dual} Banach $\mc A$-bimodule
is inner.  For example, commutative, unital C$^*$-algebras
(that is, $C(K)$ spaces with pointwise product) are amenable;
group algebras $L^1(G)$ are amenable if and only if the locally
compact group $G$ is amenable (see \cite{RundeBook} for these and
further results).

Let $\mc A$ be a Banach algebra, and turn $\mc A\proten\mc A$ into
a Banach $\mc A$-bimodule by
\[ a \cdot (b\otimes c) = ab \otimes c, \quad
(b\otimes c)\cdot a = b \otimes ca
\qquad (a\in\mc A, b\otimes c\in \mc A \proten \mc A), \]
and linearity and continuity.
Define the product map $\Delta_{\mc A} : \mc A\proten\mc A
\rightarrow \mc A$ by $\Delta_{\mc A}(a\otimes b) = ab$.
The following result, due to Johnson, can be found in \cite[Section~2.2]{RundeBook}.

\begin{theorem}\label{diag_char}
Let $\mc A$ be a Banach algebra.  Then $\mc A$ is contractible
if and only if there exists a \emph{diagonal} $\tau \in \mc A
\proten\mc A$; that is, $a\cdot\tau=\tau\cdot a$ and
$\Delta_{\mc A}(\tau)a=a$ for each $a\in\mc A$.

Similarly, $\mc A$ is amenable if and only if, for some $C>0$,
for each $\epsilon>0$ and $a_1,\ldots,a_n\in\mc A$, there exists
$\tau\in\mc A\proten\mc A$ such that $\| \tau \|_\pi \leq C$,
$\| a_i\cdot\tau - \tau\cdot a_i \|_\pi < \epsilon$ and
$\| \Delta_{\mc A}(\tau)a_i - a_i \|<\epsilon$ for
$1\leq i\leq n$.
\end{theorem}

We shall say that $\mc A$ is \emph{$C$-amenable}
if the above holds for the constant $C>0$.

We define a Banach algebra $\mc A$ to be \emph{ultra-amenable}
if every ultrapower of $\mc A$ is amenable.  It would be more
natural, in light of terms like ``super-reflexive'', to call
this super-amenable, but this term is already used by Runde in \cite{RundeBook}
for contractible algebras (as ``contractible'' has multiple
meanings as well!)  We see immediately that unital, commutative
C$^*$-algebras are ultra-amenable.

\begin{proposition}\label{surj_inherit}
Let $\mc A$ be a Banach algebra, and let $\mc I$ be a closed ideal
in $\mc A$.  When $\mc A$ is ultra-amenable, $\mc A/\mc I$
is ultra-amenable.  If $\mc I$ and $\mc A/\mc I$ are ultra-amenable,
then so is $\mc A$.  Furthermore, when $\mc A$ is ultra-amenable,
$\mc I$ is ultra-amenable if and only if $\mc I$ has a bounded
approximate identity.
\end{proposition}
\begin{proof}
Let $\mc U$ be an ultrafilter.  Then the quotient map $\mc A
\rightarrow \mc A/\mc I$ induces a quotient map $(\mc A)_{\mc U}
\rightarrow (\mc A/\mc I)_{\mc U}$, and so $(\mc A/\mc I)_{\mc U}$
is amenable by \cite[Corollary~2.3.2]{RundeBook}.  Indeed, as in the
Banach space case (compare \cite[Proposition~6.5]{Hein}), we can identify
$(\mc A/\mc I)_{\mc U}$ with $(\mc A)_{\mc U} / (\mc I)_{\mc U}$.
Hence if $(\mc A/\mc I)_{\mc U}$ and $(\mc I)_{\mc U}$ are
both amenable, then so is $(\mc A)_{\mc U}$ by \cite[Theorem~2.3.10]{RundeBook}.

As $(\mc A)_{\mc U}$ is amenable and $(\mc I)_{\mc U}$ is an ideal
in $(\mc A)_{\mc U}$, by \cite[Theorem~2.3.7]{RundeBook}, $(\mc I)_{\mc U}$
is amenable if and only if it has a bounded approximate identity.
By Proposition~\ref{bai_ultra}, this is equivalent to $\mc I$ having a bounded
approximate identity, as required.
\end{proof}

\subsection{Diagonal-like constructions}

Instead of working with the definition of amenability, it is
common to work with \emph{approximate} or \emph{virtual} diagonals;
see \cite[Chapter~2]{RundeBook} (and Theorem~\ref{diag_char} above).
In this section, we provide a similar characterisation of ultra-amenability.

\begin{definition}
Let $\mc A$ be a Banach algebra, and let $n\geq 1$.  Let
$S_n(\mc A)$ be the collection of subsets of size $n$ of the unit sphere
of $\mc A$.

Let $C>0$, $\epsilon>0$ and $n\geq 1$.  For $A\subseteq S_n(\mc A)$,
we say that $A\in D_n(\mc A,C,\epsilon)$ when there exists a sequence
of positive reals $(\alpha_k)$ with $\sum_k\alpha_k\leq C$, and such
that for each $F\in A$, we have that there exists
\[ \tau = \sum_{k=1}^\infty a_k \otimes b_k \in\mc A\proten\mc A \]
with
\[ \| a\cdot\tau - \tau\cdot a\|_\pi \leq\epsilon,\quad
\| \Delta_{\mc A}(\tau) a-a \|\leq\epsilon \qquad (a\in F), \]
and with $\|a_k\| \|b_k\| \leq \alpha_k$ for each $k$.
\end{definition}

\begin{proposition}\label{first_ultra_amen}
Let $\mc A$ be a Banach algebra, let $\mc U$ be an ultrafilter on an
index set $I$, and let $C>0$ be a constant.
Then the following are equivalent:
\begin{enumerate}
\item $(\mc A)_{\mc U}$ is $C$-amenable;
\item for each $n\geq 1$, each $\epsilon>0$, and each map
$\gamma:I\rightarrow S_n(\mc A)$, there exists a sequence of
positive reals $(\alpha_j)$, with $\sum_j \alpha_j \leq C$, and there
exists $B\in\mc U$ such that for each $i\in B$, there exists
$\tau = \sum_{j=1}^\infty b_j\otimes c_j \in \mc A\proten\mc A$ with
\[ \|a\cdot\tau - \tau\cdot a\|_\pi < \epsilon, \quad
\|\Delta_{\mc A}(\tau) a - a\| < \epsilon \qquad (a\in\gamma(i)), \]
and with $\|b_j\|\|c_j\|\leq\alpha_j$ for $j\geq 1$.
\item for each $n\geq 1$, each $\epsilon>0$, and each map
$\gamma:I\rightarrow S_n(\mc A)$, there exists $A\in D_n(\mc A,C,\epsilon)$
with $\gamma^{-1}(A) \in \mc U$.
\end{enumerate}
\end{proposition}
\begin{proof}
By definition, $(\mc A)_{\mc U}$ is amenable if and only
if there exists $C>0$ such that, for each $\epsilon>0$
and each $a^{(1)},\cdots,a^{(n)} \in (\mc A)_{\mc U}$,
there exists $\tau \in (\mc A)_{\mc U} \proten (\mc A)_{\mc U}$
such that
\[ \| a^{(k)}\cdot\tau - \tau\cdot a^{(k)} \|_\pi < \epsilon, \quad
\| \Delta_{(\mc A)_{\mc U}}(\tau) a^{(k)} - a^{(k)} \| < \epsilon
\qquad (1\leq k\leq n), \]
and with $\| \tau \|_\pi \leq C$.
Now, we may suppose that the $a^{(k)}$ are distinct, and,
by a perturbation argument, that $a^{(j)}_i \not= a^{(k)}_i$ for
each $i\in I$ and $j\not=k$, while $\| a^{(k)}_i \| = \| a^{(k)} \|$
for each $i\in I$.  As we are free to vary $\epsilon>0$,
it is enough to consider the case when $\| a^{(k)} \|=1$ for each
$k$.  Thus the choice of the family $\{ a^{(k)} : 1\leq k\leq n \}$
corresponds to a choice of a map $\gamma:I\rightarrow S_n(\mc A)$,
together with some ordering.

As explained before, $\psi_0 : (\mc A)_{\mc U} \proten
(\mc A)_{\mc U} \rightarrow (\mc A \proten \mc A)_{\mc U}$
is both an $\mc A$-bimodule homomorphism, and an
$(\mc A)_{\mc U}$-bimodule homomorphism.  Furthermore, the
following diagram commutes
\[ \spreaddiagramrows{2ex} \spreaddiagramcolumns{5ex}
\xymatrix{ (\mc A)_{\mc U} \proten (\mc A)_{\mc U}
\ar[r]^-{\Delta_{(\mc A)_{\mc U}}}
\ar[d]^{\psi_0} & (\mc A)_{\mc U} \\
(\mc A \proten \mc A)_{\mc U}
\ar[ru]^-{(\Delta_{\mc A})} } \]

Thus let $\psi_0(\tau) =
(\tau_i)_{i\in I}$, so that our conditions upon $\tau$ imply that
\[ \lim_{i\rightarrow\mc U} \| a^{(k)}_i \cdot \tau_i -
\tau_i \cdot a^{(k)}_i \|_\pi < \epsilon, \quad
\lim_{i\rightarrow\mc U} \| \Delta_{\mc A}(\tau_i) a^{(k)}_i - a^{(k)}_i \|
< \epsilon \qquad (1\leq k\leq n). \]
By definition, this is so if and only if there exists $B\in\mc U$ with
\[ \| a \cdot \tau_i - \tau_i \cdot a \|_\pi < \epsilon, \quad
\| \Delta_{\mc A}(\tau_i) a - a \| < \epsilon
\qquad (i\in B, a\in\gamma(i)). \]
With reference to Proposition~\ref{image_psi_zero},
we have that (1) implies (2).

Conversely, we simply apply Proposition~\ref{image_psi_zero}
to build $\tau \in (\mc A)_{\mc U} \proten (\mc A)_{\mc U}$
out of the family $(\tau_i)$, where the condition upon the
$(\tau_i)$ implies that we can work in the image of $\psi_0$.

Finally, it is easy to see that (2) and (3) are equivalent.
\end{proof}

\begin{corollary}
Let $\mc A$ be a Banach algebra.  If $\mc A$ is contractible,
then $\mc A$ is ultra-amenable.  If $\mc A$ is ultra-amenable,
$\mc A$ is amenable.
\end{corollary}

By considering commutative C$^*$-algebras, we see that being
contractible is strictly stronger than being ultra-amenable.
We shall see below that being ultra-amenable is also strictly
stronger than being amenable.

\begin{theorem}\label{ua_char}
Let $\mc A$ be a Banach algebra.  The following are equivalent:
\begin{enumerate}
\item\label{uac:one} $\mc A$ is ultra-amenable;
\item\label{uac:two} there exists a constant $C>0$ such that for
every $n\geq1$ and $\epsilon>0$, there exists
a finite partition $S_n(\mc A) = A_1 \cup\cdots\cup A_k$ with
each $A_i\in D_n(\mc A,C,\epsilon)$.
\end{enumerate}
\end{theorem}
\begin{proof}
Suppose that (\ref{uac:two}) holds.  We shall verify condition (3) of
Proposition~\ref{first_ultra_amen}, which will show that (\ref{uac:one}) holds.
If $\mc U$ is an ultrafilter on an index set $I$ and
$\gamma:I\rightarrow S_n(\mc A)$ is a map,
then it is easily checked that
\[ \gamma^*(\mc U) := \{ A\subseteq S_n(\mc A) : \gamma^{-1}(A)\in\mc U \} \]
is an ultrafilter on $S_n(\mc A)$.  Thus there exists some $k$ with
$A_k\in\gamma^*(\mc U)$, that is, $\gamma^{-1}(A_k)\in\mc U$, as
required.

Suppose that (\ref{uac:one}) holds.  We first introduce a little notation.
For an ultrafilter $\mc U$ on an index set $I$, for a property $p$
of elements of $I$, we write
\[ \forall_{\mc U}\,i,\, p(i) \quad \Leftrightarrow \quad
\big\{ i\in I : p(i) \text{ holds } \big\} \in \mc U. \]

For each $n\geq 1$, let $\mc U_n$ be an ultrafilter on $S_n(\mc A)$.
For each $n$, define an arbitrary injection $\iota_n:S_n(\mc A)\rightarrow
S_{n+1}(\mc A)$.  
Let $I$ be the collection of sequences $(F_n)$ where $F_n\in S_n(\mc A)$
for each $n$, and for some $N>0$ (depending on the sequence), we have that
$F_{n+1} = \iota_n(F_n)$ for $n\geq N$.  Loosely speaking, $I$ is the
collection of eventually constant sequences in $\bigcup_n S_n(\mc A)$.
Let $\mc U$ be a non-principal ultrafilter on $\mc N$.  We define an
ultrafilter $\mc V$ on $I$ by setting $K\in\mc V$ if and only if
\[ \forall_{\mc U}\,n,\, \forall_{\mc U_1}\,F_1,\, \cdots
\forall_{\mc U_n}\,F_n,\, (F_1,F_2,\cdots,F_n,F_n,\cdots)\in K. \]
It is an easy check that $\mc V$ is an ultrafilter.  Suppose that
$(\mc A)_{\mc V}$ is $C$-amenable for some constant $C>0$.

Let $n\geq 1$, and define a map $\gamma:I\rightarrow S_n(\mc A)$ by
$\gamma((F_k)) = F_n$ for $(F_k)\in I$.  Hence by Proposition~\ref{first_ultra_amen},
for each $\epsilon>0$, there exists $A\in D_n(\mc A,C,\epsilon)$ with
$\gamma^{-1}(A) \in \mc V$.  By the definition of $\mc V$, this means that
\[ \forall_{\mc U}\,m,\, \forall_{\mc U_1}\,F_1,\, \cdots
\forall_{\mc U_m}\,F_m,\, \gamma(F_1,F_2,\cdots,F_m,F_m,\cdots)\in A. \]
Hence, for some $m>n$, by the definition of $\gamma$, we see that
\[ \forall_{\mc U_1}\,F_1,\, \cdots
\forall_{\mc U_m}\,F_m,\, F_n\in A, \]
that is, $\forall_{\mc U_n}\,F_n,\, F_n\in A$, which is simply the
statement that $A\in\mc U_n$.  In conclusion, for each $n\geq1$ and
$\epsilon>0$, there exists some member of $D_n(\mc A,C,\epsilon)$ in $\mc U_n$.

If condition (\ref{uac:two}) does not hold, then
for some $n\geq1$ and $\epsilon>0$, we have that there is no
finite cover of $S_n(\mc A)$ by members of $D_n(\mc A,C,\epsilon)$.
In particular, if $\mc F = \{ S_n(\mc A) \setminus A : A\in D_n(\mc A,C,\epsilon) \}$
then, as no finite intersection of members of $\mc F$ is empty,
there exists an ultrafilter $\mc U_n$ containing $\mc F$.
However, by the previous paragraph, we know that $\mc U_n \cap
D_n(\mc A,C,\epsilon)$ is non-empty, which gives a contradiction,
as required.
\end{proof}

\subsection{C$^*$-algebras}

By throwing a lot of machinery at the problem, we can rather
easily settle the question of when a C$^*$-algebra is ultra-amenable.
In \cite{LLW}, there is a throwaway comment in the proof of Theorem~2.5
that, for a C$^*$-algebra $\mc A$, every ultrapower $(\mc A)_{\mc U}$
is amenable if and only if $\ell^\infty(\mc A,I)$ is amenable for all
index sets $I$.  We do not see why this is ``obvious'' however, as
in general $(\mc A)_{\mc U}$ is much smaller than $\ell^\infty(\mc A,I)$.
The following proof avoids this issue, and uses no more machinery than
\cite{LLW} does.

\begin{theorem}\label{cstar_case}
Let $\mc A$ be a C$^*$-algebra.  Then the following are equivalent:
\begin{enumerate}
\item\label{cstar:ua} $\mc A$ is ultra-amenable;
\item\label{cstar:bdam} $\mc A''$ is amenable;
\item\label{cstar:subh} $\mc A$ is the finite-direct sum of algebras of
   the form $C_0(K) \otimes \mathbb M_n$ for some integer $n$ and some
   locally compact Hausdorff space $K$.
\end{enumerate}
\end{theorem}
\begin{proof}
By \cite[Corollary~6.4.28]{RundeBook} and \cite[Corollary~6.4.29]{RundeBook},
we have that when a C$^*$-algebra is amenable, it has the (metric)
approximation property.  Hence, if (\ref{cstar:ua}) holds, then
for any ultrafilter $\mc U$, we have that $(\mc A)_{\mc U}$ has the
approximation property.  For a suitable choice of $\mc U$, we have that
$\mc A''$ is isometric to a complemented subspace of $(\mc A)_{\mc U}$
(see \cite[Proposition~6.7]{Hein}).  As the approximation property clearly
drops to complemented subspaces, we conclude that $\mc A''$ has the approximation
property.  By \cite[Theorem~6.1.7]{RundeBook} and \cite[Remark~6.1.9]{RundeBook},
this implies that $\mc A''$ is amenable, giving (\ref{cstar:bdam}).

When (\ref{cstar:bdam}) holds, by \cite[Theorem~6.1.7]{RundeBook}, we know
that $\mc A''$ has the form specified in (\ref{cstar:subh}), but with
$K$ compact, and such that $C_0(K) = C(K)$ is a dual space (so that $K$
is actually a \emph{hyperstonian} space).  Suppose that $\mc A''$ were
isomorphic to $C(K)$ for some hyperstonian compact space $X$.  Then $\mc A$
is commutative, and so is isomorphic to $C_0(L)$, for some locally compact
space $L$.  Then note that the bidual of $C_0(L) \otimes \mathbb M_n$
is isomorphic to $C(K) \otimes \mathbb M_n$, which is isomorphic to
$\mc A''$.  It is now clear that $\mc A$ must be isomorphic to
$C_0(L) \otimes \mathbb M_n$, showing (\ref{cstar:subh}).

Finally, when (\ref{cstar:subh}) holds, it is clear that $\mc A$ is
ultra-amenable, giving (\ref{cstar:ua}).
\end{proof}

It seems unlikely that (1) and (2) are equivalent for a general
Banach algebra $\mc A$, but we do not currently have a counter-example.
It would, of course, be nice to be able to prove the above without
using so much machinery, even for certain classes of C$^*$-algebra.

\subsection{Group algebras}

Let $G$ be a locally compact group, and form the Banach algebra $L^1(G)$.
See \cite[Section~3.3]{Dales} or \cite[Section~1.9]{Palmer1}, for example,
for details about this class of algebras.  We shall make use of the
concept of the \emph{almost periodic}, or \emph{Bohr}, compactification
for a group $G$.  However, it makes sense for us to develop these ideas
first for general Banach algebras.

Above, for a Banach algebra $\mc A$, we defined the space of
weakly almost periodic functionals of $\mc A$, denoted by $\wap(\mc A')$.
If we insist that the map $L_\mu$, for $\mu\in\mc A'$,
\[ L_\mu: \mc A\rightarrow\mc A'; \quad a\mapsto a\cdot\mu
\qquad (a\in\mc A), \]
is compact, and not just weakly-compact, we arrive at the definition
of an \emph{almost periodic functional}, denoted by $\mu\in\ap(\mc A')$.
Clearly $\ap(\mc A')\subseteq\wap(\mc A')$, and it is easy to show that
$\ap(\mc A')$ is a closed submodule of $\mc A'$.

Let $\mc A$ be a Banach algebra, and let $\mc U$ be an ultrafilter.
Then we can define
\[ \sigma_{{\ap}}:(\mc A)_{\mc U}\rightarrow\ap(\mc A')',
\quad \ip{\sigma_{{\ap}}((a_i))}{\mu} = \lim_{i\rightarrow\mc U} \ip{\mu}{a_i}
\qquad (\mu\in\ap(\mc A'), (a_i)\in(\mc A)_{\mc U}). \]
It is clear that $\sigma_{{\ap}}$ is norm-decreasing, and for suitable
$\mc U$, $\sigma_{{\ap}}$ is a surjection (compare \cite[Proposition~6.7]{Hein}).
As remarked upon in \cite[Section~5]{Daws}, $\sigma_{\ap}$ is easily seen to
be an algebra homomorphism.

\begin{proposition}
Let $\mc A$ be an ultra-amenable Banach algebra.  Then
$\ap(\mc A')'$ is amenable.
\end{proposition}
\begin{proof}
This is immediate, as $\ap(\mc A')'$ can be identified with
a quotient of a suitable ultrapower of $\mc A$.
\end{proof}

We note that for some algebras $\mc A$, $\ap(\mc A')$ can be
trivial.  For example, let $E$ be an infinite dimensional
Banach space, and let $\mc A=\mc A(E)$ be the algebra of
approximable operators on $E$.  Then the dual of $\mc A(E)$ is
$\mc I(E')$, then space of integral operators on $E$.  For
$U\in\mc I(E')$, $T\in\mc A(E)$ and $\mu\otimes x\in\mc A(E)$, we see that
\[ \ip{(\mu\otimes x)\cdot U}{T} = \ip{U}{T(\mu\otimes x)}
= \ip{U}{\mu\otimes T(x)} = \ip{U(\mu)}{T(x)}
= \ip{U(\mu)\otimes x}{T}. \]
Let $U\in\mc I(E')$ be non-zero, and choose $\mu\in E'$ with
$\|\mu\|=1$ and such that $\lambda=U(\mu)$ is non-zero.
Define a map $R_\lambda:E\rightarrow
\mc I(E')$ by $R_\lambda(x) = \lambda\otimes x$, so that $R_\lambda$
is an isomorphism onto its range.  We hence have that
$(\mu\otimes x)\cdot U = R_\lambda(x)$ for each $x\in E$.
As
\[ \big\{ S\cdot U : S\in\mc A(E), \|S\|\leq 1 \big\}
\supseteq \big\{ U(\mu)\otimes x : x\in E, \|x\|\leq 1 \big\}
= \big\{ R_\lambda(x) : \|x\|\leq 1 \big\}, \]
we see if $U\in\ap(\mc A')$, then $R_\lambda$ is compact,
which implies that $E$ is finite-dimensional, a contradiction.
We conclude that $\ap(\mc A')=\{0\}$.

When $\mc A=L^1(G)$, however, $\ap(\mc A')$ is often large.
We write $\ap(G)$ for $\ap(\mc A')$.
For $g\in G$ define
\[ L_g:L^\infty(G)\rightarrow L^\infty(G); \quad
L_g(\mu)(h) = \mu(g^{-1}h) \qquad (\mu\in L^\infty(G), h\in G). \]
When $G$ is discrete, we see that for $\mu\in L^\infty$, the set
$\{ a\cdot\mu : a\in\ell^1(G), \|a\|\leq 1 \}$ is contained in
the closure of the absolutely convex hull of the set
$\{ L_g(\mu) : g\in G \}$.  Hence $\mu\in\ap(G)$ if and only
if the set of translates $\{ L_g(\mu) : g\in G \}$ is
relatively compact in $L^\infty(G)$.  By a more intricate argument,
we can show that this is true for general $G$ (compare with the
argument in \cite{Ulger}).

For more details on $\ap(G)$, see for example \cite[Section~3.2.16]{Palmer1}
and \cite[Theorem~12.4.15]{Palmer2}.  We have that $\ap(G)$ is a unital
C$^*$-subalgebra of $L^\infty(G)$, so that $\ap(G) = C(K)$, for some
compact space $K$.  We denote $K$ by $G^{\ap}$, so that
$\ap(\mc A')' = M(G^{\ap})$, the space of measures on $G^{\ap}$.  Each
member of $G$ induces a character on $\ap(G)$, and this leads to a canonical
map $G\rightarrow G^{\ap}$ which has dense range.  We can use this to extend
the product on $G$ to a product on $G^{\ap}$, which turns $G^{\ap}$ into a
compact group.  We can check that the induced product on $\ap(\mc A')'$
agrees with the convolution product on $M(G^{\ap})$.  In fact, $G^{\ap}$
has the following universal property.  If $H$ is any compact group, and
$\phi:G\rightarrow H$ is a continuous homomorphism, then $\phi$ factors
through to canonical map $G\rightarrow G^{\ap}$.  In this sense, $G^{\ap}$
is the maximal compact group which contains a dense homomorphic image of $G$.

In general, the canonical map $G\rightarrow G^{\ap}$ need not be
an injection.  When it is, equivalently, when $\ap(G)$ separates
the points of $G$, we say that $G$ is \emph{maximally almost
periodic}.  Obviously all compact groups are maximally almost
periodic, and abelian groups $G$ are also maximally almost periodic,
which follows as the characters of $G$ separate the points of $G$.
At the other extreme, if $\ap(G)$ is the linear span of the constant
function in $L^\infty(G)$, then $G$ is \emph{minimally almost periodic}.
This is equivalent to the statement that the only continuous
homomorphisms from $G$ to a compact group are trivial.

There exist amenable, minimally almost periodic groups, as the
following example, due to George Willis, shows.  Let
$\operatorname{FSym}(\mathbb N)$ be the collection of permutations of
the natural numbers which fix all but finitely many elements.  Let
$G=\Alt(\mathbb N)$ be the index 2 subgroup of even permutations.
The $\Alt(\mathbb N)$ is a simple group (see \cite[Corollary~3.3A]{DM})
and as it is the direct limit of finite groups, it is amenable.
Suppose that $\ap(G)$ is not trivial, so that $G$ admits some
non-trivial homomorphism into a compact group.  By the representation
theory of compact groups (essentially, the Peter-Weil theorem),
it follows that there is a non-trivial homomorphism of $G$ into
a matrix group $GL_n(\mathbb C)$.  As $G$ is simple, such a
homomorphism is injective.  By a theorem of Tits (see \cite[Theorem~3.10]{Pat})
it follows that $G$ contains a normal, solvable group of finite index.
This is a contraction, and so we see that $\ap(G)$ is trivial.

\begin{theorem}\label{ab_com_case}
Let $G$ an infinite compact group, or an infinite abelian locally
compact group.  Then $L^1(G)$ is amenable, but not ultra-amenable.
\end{theorem}
\begin{proof}
It is well know that for such groups $G$, $L^1(G)$ is amenable.
Suppose that $L^1(G)$ is ultra-amenable, so by the above proposition,
we see that $M(G^{\ap})$ is amenable.  By \cite{DGH}, this implies
that $G^{\ap}$ is amenable and discrete, but as $G^{\ap}$ is
compact, we have that $G^{\ap}$ is finite.  However, we remarked
above that $G$ is maximally periodic, and so we see that $G$ is
also finite, a contradiction.
\end{proof}

\subsection{Discrete Group algebras}

We wish to develop a little theory for \emph{discrete} group
algebras.  This will motivate some technical and obtuse constructions
below, which will work for many non-compact group.  Let
$G$ be a discrete group, and consider the Banach algebra $\ell^1(G)$.
We write $\delta_g$ for the point mass at $g\in G$, so every
$a\in\ell^1(G)$ can be written as $a=\sum_{g\in G} a_g \delta_g$ for
some family of scalars $(a_g)_{g\in G}$ such that $\|a\|=\sum_g |a_g|$.
Let $\mc U$ be an ultrafilter on an index set $I$, so we can form the
ultrapower $(\ell^1(G))_{\mc U}$.  We can also form the \emph{ultrapower}
of $G$, denoted by $(G)_{\mc U}$.  This is the set of all families
$(g_i)_{i\in I}$ of elements of $G$, quotiented by the
equivalence relation
\[ (g_i) \sim (h_i) \quad\Leftrightarrow\quad
\{ i\in I : g_i=h_i \}\in\mc U. \]
Then $(G)_{\mc U}$ becomes a (discrete) group for the pointwise
product, and we have a canonical map $G\rightarrow (G)_{\mc U}$
formed by sending $g\in G$ to the constant family $(g)$.

For $1\leq p<\infty$, define a map $\psi_p:\ell^p((G)_{\mc U})
\rightarrow (\ell^p(G))_{\mc U}$ by
\[ \psi_p( \delta_g ) = ( \delta_{g_i} )_{i\in I}
\qquad (g = (g_i) \in (G)_{\mc U}). \]
If $(g_i) \sim (h_i)$ in $(G)_{\mc U}$ then
\[ \{ i\in I : \|\delta_{g_i} - \delta_{g_i}\|=0 \}
= \{ i\in I : g_i=h_i \}\in\mc U, \]
showing that $\psi_p$ is well-defined.  An analogous calculation
shows that if we extend $\psi_p$ by linearity and continuity,
then $\psi_p$ is an isometry onto its range.  Let $\psi_0:c_0((G)_{\mc U})
\rightarrow (c_0(G))_{\mc U}$ be the analogous map.  Then it is easy
to check that when $1<p<\infty$ and $p^{-1} + p'^{-1}=1$, then
$\psi_{p'}' \circ \psi_p$ is the identity on $\ell^p((G)_{\mc U})$.
Similarly $\psi_0' \circ \psi_1$ is the identity on $\ell^1((G)_{\mc U})$,
where as usual, we treat $(\ell^1(G))_{\mc U}$ as a subspace of
$(c_0(G))_{\mc U}'$.  For $1\leq p<\infty$, we can hence identify
$\ell^p((G)_{\mc U})$ with a $1$-complemented subspace of $(\ell^p(G))_{\mc U}$.
This identification respects the identification of $\ell^p(G)$ in either
$\ell^p((G)_{\mc U})$ or $(\ell^p(G))_{\mc U}$.
Let $I_p$ be the obvious complementary subspace to $\ell^p((G)_{\mc U})$ in
$(\ell^p(G))_{\mc U}$, that is, $I_p$ is the kernel of $\psi_{p'}'$, or
$\psi_0'$, as appropriate.

\begin{lemma}
We may identify $I_p$ with the collection of equivalence classes
in $(\ell^p(G))_{\mc U}$ represented by sequences $(x_i)_{i\in\mathbb N}$
with $\lim_{i\rightarrow\mc U} \|x_i\|_\infty=0$.  Furthermore, $I_1$ is
an ideal in the algebra $(\ell^1(G))_{\mc U}$.
\end{lemma}
\begin{proof}
Suppose that $(x_i)\in (\ell^p(G))_{\mc U}$ is such that $\lim_{i\rightarrow\mc U}
\|x_i\|_\infty>0$.  For each $i$, let $x_i = (x^{(i)}_g)_{g\in G}
\in \ell^p(G)$.  Hence there exists $\delta>0$ and $U\in\mc U$
with a function $k:U\rightarrow G$ such that $|x^{(i)}_{k(i)}|\geq\delta$
for $i\in U$.  Extend $k$ to $I$ in an arbitrary way, so we see that
\[ \lim_{i\rightarrow\mc U} | \ip{\delta_{k(i)}}{x_i} | \geq\delta, \]
and so $(x_i)$ does not annihilate $(\delta_{k(i)}) \in \ell^p((G)_{\mc U})$.
Hence $(x_i)\not\in I_p$, as required.

Now suppose that $(x_i)\in (\ell^p(G))_{\mc U}$ is such that
$\lim_{i\rightarrow\mc U} \|x_i\|_\infty = 0$.  For any map $k:I\rightarrow G$,
we see that
\[ \lim_{i\rightarrow\mc U} | \ip{\delta_{k(i)}}{x_i} |
\leq \lim_{i\rightarrow\mc U} \|x_i\|_\infty = 0, \]
so that $(x_i)$ annihilates $(\delta_{k(i)}) \in \ell^p((G)_{\mc U})$.
By linearity and continuity, $(x_i)\in I_p$, as required.

Finally, consider $(x_i)\in I_1$ and let $(a_i)\in(\ell^1(G))_{\mc U}$, so that
\begin{align*} \lim_{i\rightarrow\mc U} \|a_ix_i\|_\infty &= \lim_{i\rightarrow\mc U}
\sup_{g\in G} \Big| \sum_{h\in G} a^{(i)}_h x^{(i)}_{h^{-1}g} \Big|
\leq \lim_{i\rightarrow\mc U}
\sup_{g\in G} \sum_{h\in G} |a^{(i)}_h| |x^{(i)}_{h^{-1}g}| \\
&\leq \lim_{i\rightarrow\mc U} \|x_i\|_\infty \|a_i\|_1 = 0,
\end{align*}
so that $(a_ix_i)\in I_1$.  Hence $I_1$ is a left-ideal, and
similarly $I_1$ is a right-ideal.
\end{proof}

\begin{theorem}\label{dis_case}
Let $G$ be an infinite discrete group, and let $\mc U$ be a
countably incomplete ultrafilter.  Then $(\ell^1(G))_{\mc U}$
is not amenable.  In particular, $\ell^1(G)$ is not ultra-amenable.
\end{theorem}
\begin{proof}
As $\mc U$ is countably incomplete, we may suppose that $\mc U$ is an
ultrafilter on $\mathbb N$ (compare with the proofs of Theorem~6.3 or
Proposition~7.1 in \cite{Hein}).  Suppose that $(\ell^1(G))_{\mc U}$ is
amenable.  Then $I_1$ is a complemented ideal in an amenable Banach
algebra, and so $I_1$ is amenable (see \cite[Theorem~2.3.7]{RundeBook}).
In particular, $I_1$ has a bounded approximate identity.  We shall
show that this leads to contradiction, as required.

Let $H$ be a countably infinite subgroup of $G$, and choose a sequence
$(g_n)$ in $H$ as follows.  Let $g_1$ be arbitrary.  Suppose we have
chosen $g_1,\cdots,g_n$.  Consider the set
\[ B_n = (g_i)_{i=1}^n \cup (g_ig_j^{-1}g_k)_{1\leq i,j,k\leq n}, \]
which is finite.  We simply choose $g_{n+1}\in H\setminus B_n$.

Let $e_G$ be the unit of $G$, and let $g\in H$ with $g\not=e_G$.  Suppose
that $g_ig=g_j$ for some $i,j\geq 1$.  Let $t\geq 1$ be minimal such
that, for some $m\geq 1$, we have that $g=g_t^{-1}g_m$ or $g=g_m^{-1}g_t$.
By minimality, $t<m$, and so $g_ig\in B_m$ for $1\leq i\leq m$.  Thus
$g_n\not=g_ig$, equivalently $g\not=g_i^{-1}g_n$, for any $n>m$ and
$i\leq m$.  Similarly, $g_ig^{-1}\in B_m$ for $i\leq m$, and so
$g \not= g_n^{-1}g_i$ for $n>m$ and $i\leq m$.

Suppose that $g=g_t^{-1}g_m$, and let $r,s\geq 1$ be such that $g=g_r^{-1}g_s$
and $(r,s)\not=(t,m)$.  If $r=t$ then $g=g_t^{-1}g_m = g_t^{-1}g_s$ so that
$g_m=g_s$, that is, $m=s$, a contradiction.  Similarly $m\not=s$.
By minimality, $r,s\geq t$, so actually $r>t$.  By the above, $r,s\leq m$,
so actually $s<m$.  If $r<m$ then $g_tg\in B_{m-1}$ and so $g_m\not=
g_tg$, a contradiction.  Hence $r=m$ and so $g=g_t^{-1}g_m=g_m^{-1}g_s$,
so that $g_s = g_m g_t^{-1} g_m$, and hence $s$ must be unique.
A similar argument works when $g=g_m^{-1}g_t$.  We conclude that 
\[ \big| \big\{ k\geq 1 : g_k g = g_m \text{ for some }m\geq 1 \big\}\big|
\leq 2 \qquad (g\in H, g\not=e_G). \]

Set $x = \chi_{\{g_n : n\in\mathbb N\}} \in \ell^\infty(G)$, so that $\|x\|=1$.
For each $n\in\mathbb N$, let
\[ a_n = \frac1n \sum_{i=1}^n \delta_{g_i} \in \ell^1(G), \]
so that $\|a_n\|_1 = 1$.  We have that
\[ \ip{x\cdot a_n}{\delta_{e_G}} = \frac1n\sum_{i=1}^n \ip{x}{\delta_{g_i}} = 1. \]
Now let $g\in H$ with $g\not=e_G$.  Then, by the above,
\[ |\ip{x\cdot a_n}{\delta_g}| = \frac1n\sum_{i=1}^n \ip{x}{\delta_{g_ig}}
=\frac1n\big|\big\{ 1\leq i\leq n : g_ig=g_k \text{ for some }k\geq1 \big\}\big|
\leq \frac2n. \]
If $g\in G\setminus H$, then clearly $\ip{x}{\delta_{g_ig}}=0$, as
$g_ig\not\in H$, for any $i\geq 1$.  Let $a=(a_n)\in (\ell^1(G))_{\mc U}$,
so clearly $a\in I_1$, and we see that $x\cdot a = e_{e_G}$ in $(\ell^\infty(G))_{\mc U}$.
Here $e_{e_G}$ refers to the point mass at $e_G$, namely the same function as
$\delta_{e_G}$, but now treated as a member of $c_0(G) \subseteq \ell^\infty(G)$.

By assumption, $I_1$ has a bounded approximate identity, so in
particular, there exists $b=(b_n)\in I_1$ with $\|a-ab\|<1/2$.
Hence $|\ip{x\cdot a}{\delta_{e_G}-b}| = |\ip{x}{a-ab}| < 1/2$, as
$\|x\|_\infty=1$.  However, from the above, $\ip{x\cdot a}{\delta_{e_G}-b}
= \ip{e_{e_G}}{\delta_{e_G} - b} = 1-\lim_{i\rightarrow\mc U} \ip{e_{e_G}}{b_i} = 1$,
as $b\in I_1$.  This contradiction completes the proof.
\end{proof}

\subsection{General groups}

We start by making some observations about quotients of groups.

\begin{proposition}\label{quot_group}
Let $G$ be a locally compact group such that $L^1(G)$ is ultra-amenable,
and let $H$ be a closed normal subgroup of $G$.  If $G/H$ is compact,
abelian or discrete, then $G/H$ is finite.
\end{proposition}
\begin{proof}
As detailed in \cite[Section~1.9.12]{Palmer1}, we have a surjective
algebra homomorphism $L^1(G)\rightarrow L^1(G/H)$.  By
Proposition~\ref{surj_inherit}, we see that $L^1(G/H)$ is ultra-amenable.
The result now follows from Theorem~\ref{ab_com_case} and
Theorem~\ref{dis_case}.
\end{proof}

In particular, by considering the modular function of the Haar
measure on $G$ (see \cite[Section~1.9]{Palmer1} or
\cite[Section~3.3]{Dales}), we see that if $L^1(G)$ is ultra-amenable,
then $G$ is unimodular, as otherwise, a quotient of $G$ would be
isomorphic to an infinite subgroup of $(\mathbb R_{>0},\times)$,
which is abelian.  Similarly, if $L^1(G)$ is ultra-amenable,
then the derived subgroup $G'$ of $G$ (see \cite[Section~12.1]{Palmer2}) must
be ``large'', in the sense that $G/G'$ is finite.

We now wish to generalise the arguments used above for discrete
groups.  Let $G$ be a locally compact group, and let $\mc U$ be
an ultrafilter on an index set $I$.  Let $I_1 \subseteq (L^1(G))_{\mc U}$
be the collection of elements $x\in (L^1(G))_{\mc U}$ such that $x$
has a representation of the form $(x_i)_{i\in I}$ where $x_i\in C_0(G)$
for each $i$, and $\lim_{i\rightarrow\mc U} \|x_i\|_\infty=0$.  It is easy
to verify that $I_1$ is a subspace of $(L^1(G))_{\mc U}$, and as
$C_0(G)\cap L^1(G)$ is dense in $L^1(G)$, it follows that $I_1$ is
closed.

\begin{lemma}
With notation as above, $I_1$ is an ideal in $(L^1(G))_{\mc U}$.
\end{lemma}
\begin{proof}
Let $x\in I_1$, so $x=(x_i)$ with $x_i$ as above.  Let $a=(a_i)\in(L^1(G))_{\mc U}$,
where by density, we may suppose that $a_i\in C_0(G)$ for each $i$.
Then $a_ix_i\in C_0(G)$ for each $i$, and we have that $\|a_ix_i\|_\infty
\leq \|a_i\|_\infty \|x_i\|_1$, from which it follows that $ax=(a_ix_i)
\in I_1$, so we see that $I_1$ is a left-ideal.  Similarly
$I_1$ is a right-ideal.
\end{proof}

We cannot, in general, show that $I_1$ is complemented.  However,
we shall show that $I_1$ is \emph{weakly-complemented}.  That is,
\[ I_1^\perp = \{ \mu\in (L^1(G))_{\mc U}' : \ip{\mu}{x}=0 \
(x\in I_1) \} \]
is complemented in $(L^1(G))_{\mc U}'$.

\begin{proposition}
With notation as above, $I_1$ is weakly-complemented.
\end{proposition}
\begin{proof}
We shall sketch this, as the details are very similar to ideas used
to deal with dual Banach algebras in Section~\ref{ultra_dual_ba} above.
Let $\phi$ be the composition of the isometric inclusions
\[ (L^1(G))_{\mc U} \rightarrow (M(G))_{\mc U} = (C_0(G)')_{\mc U}
\rightarrow (C_0(G))_{\mc U}', \]
so that $\phi' : (C_0(G))_{\mc U}'' \rightarrow (L^1(G))_{\mc U}'$
is a quotient map (or metric surjection).  As $(C_0(G))_{\mc U}$ is
a C$^*$-algebra, we see that $(C_0(G))_{\mc U}''$ is a von Neumann
algebra.  We can verify that the kernel of $\phi'$ is an ideal in
$(C_0(G))_{\mc U}''$, and so $(L^1(G))_{\mc U}'$ becomes a C$^*$-algebra.
As $(L^1(G))_{\mc U}'$ is a dual space, we see that $(L^1(G))_{\mc U}'$
is a commutative von Neumann algebra.  Notice that $(L^1(G))_{\mc U}'$
isometrically contains $(L^\infty(G))_{\mc U}$, and it is not hard to show
that $(L^\infty(G))_{\mc U}$ becomes a $*$-subalgebra of $(L^1(G))_{\mc U}'$.

We can check that $I_1^\perp$ is a $*$-subalgebra of $(L^1(G))_{\mc U}'$.
As $I_1^\perp$ is weak$^*$-closed, we see that $I_1^\perp$ is
a commutative von Neumann algebra, and is hence \emph{injective}
(see \cite[Section~6.2]{RundeBook}).  Hence there is a (contractive)
projection $(L^1(G))_{\mc U}' \rightarrow I_1^\perp$, as required.
\end{proof}

We now make a temporary definition.  Let $G$ be a non-compact,
locally compact group.  We shall say that $G$ is \emph{relatively-[IN]}
if there exists a compact, symmetric, non-null (with respect to Haar
measure) set $K$ in $G$ and a subset $A\subseteq G$ whose closure is
not compact, such that $aK=Ka$ for each $a\in A$.  We shall say that
$(K,A)$ is a \emph{witness}.  We recall that if
we can take $A=G$, then $G$ is an [IN]-group (see \cite[Section~12.1.8]{Palmer2}
for further details of this class of groups).

\begin{lemma}
Let $G$ be a relatively-[IN] group witnessed by $(K,A)$.  Let $A_0$
be the closed subgroup generated by $A$.  Then $aK=Ka$ for $a\in A_0$.
Let $K_0$ be the closed subgroup generated by $K$.  Then also
$a K_0 = K_0 a$ for $a\in A_0$.  Finally, if $G_0$ is the closed subgroup
of $G$ generated by $K_0$ and $A_0$, then $gK_0=K_0g$ for $g\in G_0$.
\end{lemma}
\begin{proof}
As $K$ is symmetric, we see that if $g,h\in A$, then $g^{-1}K = Kg^{-1}$,
and $ghK = gKh = Kgh$.  Hence $K$ is invariant under the action of
the subgroup generated by $A$.  Let $(a_\alpha)$ be a net in the 
subgroup generated by $A$ converging to $a\in G$.  For $k\in K$, the
net $(a_\alpha k a_\alpha^{-1})$ is in $K$ and converges to $aka^{-1}$,
so as $K$ is closed, $aka^{-1}\in K$, and so we see that $a K = K a$.
Thus $K$ is invariant under the action of $A_0$.

As $K$ is symmetric, the subgroup generated by $K$ is simply
$\bigcup_{n\geq 1} K^n$, and it is clear that this is invariant
under the action of $A_0$.  Let $a\in A_0$ and let $(k_\alpha)$ be
a net in the subgroup generated by $K$ tending to $k\in G$.  Then
$ak_\alpha a^{-1} \rightarrow aka^{-1}$, showing that $aka^{-1}\in K_0$.
Hence $K_0$ is invariant under the action of $A_0$.

Let $G_1$ be the subgroup of $G$ generated by $A_0$ and $K_0$.  As $K_0$
is invariant under $A_0$, we see that $G_1 = A_0 K_0 = K_0 A_0$, and
so clearly $K_0$ is invariant under the action of $G_1$.  Again, a
continuity argument shows that the same holds for $G_0$.
\end{proof}

We shall see later that we really only care about subgroups, and
so the above shows that being a relatively-[IN] group is rather
similar to having a closed subgroup which is [IN].  However, we
cannot in general show that $K$ (as opposed to $K_0$) is
invariant for $G_0$.  Indeed, $K$ is $G_0$-invariant if and only
if $K$ is $K$-invariant ($kK=Kk$ for each $k\in K$).

However, we can always assume that for a witness $(K,A)$, we
have that $A$ is a closed non-compact subgroup of $G$.

\begin{theorem}\label{in_case}
Let $G$ be a non-compact group which is relatively-[IN].  Let $\mc U$
be a countably incomplete ultrafilter.  Then $(L^1(G))_{\mc U}$
is not amenable.  In particular, $L^1(G)$ is not ultra-amenable.
\end{theorem}
\begin{proof}
As in the proof of Theorem~\ref{dis_case} above,
we may suppose that $\mc U$ is an ultrafilter on the index set $\mathbb N$.
If $(L^1(G))_{\mc U}$ is amenable, then as $I_1$ is a weakly-complemented
ideal, then $I_1$ is amenable as well, by \cite[Theorem~2.3.7]{RundeBook}.
Again, in particular, $I_1$ has a bounded approximate identity.  Let
$(K,A)$ be a witness to the fact that $G$ is relatively-[IN].

We choose a sequence $(g_n)$ in $A$ as follows.  Let $g_1\in A$
be arbitrary.  Suppose we have chosen $g_1,\cdots,g_n$, and let
\[ B_n = \bigcup_{i=1}^n g_i KK \cup \bigcup_{i=1}^n g_i K^4
\cup \bigcup_{i=1}^n g_i K^8 \cup
\bigcup_{i,j,k=1}^n g_i g_j^{-1} g_k K^4, \]
so that $B_n$ is compact in $G$.  We can hence choose $g_{n+1}\not\in B_n$,
as $A$ is not compact.  Then, for $k\leq n$, we see that $g_{n+1} \not\in
g_k KK = KK g_k$, so that $g_{n+1} g_k^{-1}\not\in KK$.
Similarly, as $(KK)^{-1}=KK$, we see that $g_kg_{n+1}^{-1}\not\in KK$ for
$k\leq n$.

Let $g\not\in KK$, and let $t\geq1$ be minimal such that, for some
$m\geq1$, we have that $g\in g_m g_t^{-1}KK$ or $g\in g_t g_m^{-1}KK$.
As $g\not\in KK$, by minimality, we have that $t<m$.
Let $r,s\geq 1$ with $g\in g_rg_s^{-1}KK$, so by minimality, $r,s\geq t$,
so actually, $r>t$ or $s>t$.

Suppose that $r,s<m$, so that either
$g\in g_mg_t^{-1}KK$, so that
$g_m \in g KK g_t \subseteq g_r g_s^{-1} K^4 g_t =
g_r g_s^{-1} g_t K^4 \subseteq B_{m-1}$, a contradiction; or
$g\in g_tg_m^{-1}KK = KKg_tg_m^{-1}$, so that
$g_m \in g^{-1} KK g_t \subseteq g_s g_r^{-1} K^4 g_t
= g_sg_r^{-1}g_t K^4\subseteq B_{m-1}$, a contradiction.
Hence $r\geq m$ or $s\geq m$.

As $r\not=s$, either $r<s$ or $s<r$.  If $r<s$, then the argument
in the previous paragraph shows that we do not have that $m,t<s$,
that is, $m\geq s$.  Similarly, if $s<r$ then $m\geq r$.  We
conclude that $m\geq\max(r,s)\geq m$, so that $m=\max(r,s)$.

Suppose that $g\in g_m g_t^{-1} KK$.  If $m=r>s$, then
$g_m^{-1}g \in g_t^{-1}KK \cap g_s^{-1}KK$.  If $t>s$ then
$g_t \in K^4 g_s$, a contradiction, so by symmetry, $t=s$.
Otherwise $m=s>r$, in which case
$g\in g_m g_t^{-1}KK \cap g_r g_m^{-1}KK$ so that $g_m g_t^{-1} g_m
\in g_r K^4$.  Suppose that there exists $r'$ with
$r'<s$ and $g\in g_{r'} g_m^{-1}KK$, so that $g_m g_t^{-1} g_m
\in g_{r'} K^4$.  Suppose that $r>r'$, so that $g_r\in g_{r'} K^8$,
a contradiction, so by symmetry, we conclude that $r=r'$.
Hence, $r$ is unique.

An analogous argument works when $g\in g_t g_m^{-1}KK$, showing
that in all cases,
\[ \big|\big\{ 1\leq l\leq n : g\in g_k g_l^{-1}KK \text{ for some }
k\geq 1 \big\}\big| \leq 2 \qquad (n\geq 1, g\in G\setminus KK). \]

For a measurable subset $B\subseteq G$ we let $\chi_B$ be the
indicator function of $B$, so that $\chi_B\in L^\infty(G)$.
When $B$ has finite measure, we have that $\chi_B\in L^1(G)$; write
$\chi^1_B$ in this case.
Let $x\in L^\infty(G)$ be defined by the following formal sum
\[ x = \sum_{n=1}^\infty \chi_{g_n K}, \]
which makes sense, as by construction, $g_n K \cap g_m K = \emptyset$ when
$n\not=m$.  We see that $\|x\|_\infty = 1$.  For each $n\geq1$, let
\[ a_n = \frac{1}{n|K|} \sum_{k=1}^n \chi^1_{g_k K}\in L^1(G), \]
so that $\|a_n\|=1$, and hence $a=(a_n)\in (L^1(G))_{\mc U}$.
Notice that $\ip{x}{a_n} = 1$ for all $n$.
Define $f\in L^\infty(G)$ by $f=|K|^{-1}\chi^1_K\cdot\chi_K$, so that
\[ f(g) = \frac{1}{|K|} \int_G \chi_K(gh) \chi^1_K(h) \ dh
= \frac{1}{|K|} \int_K \chi_{g^{-1}K}(h) \ dh =
\frac{| K \cap g^{-1}K |}{|K|} \qquad (g\in G). \]
Then $f$ is continuous, $\|f\|_\infty\leq1$ and $f(e_G)=1$.
Furthermore, $f(g)\not=0$ only when $g\in KK$.
For $s,t\in G$, we see that
\begin{align*} |K|^{-1} (\chi^1_{sK} \cdot \chi_{tK})(g) &=
\frac{1}{|K|} \int_G \chi_{tK}(gh) \chi^1_{sK}(h) \ dh
= \frac{1}{|K|} \int_G \chi_{K}(t^{-1}gh) \chi^1_{K}(s^{-1}h) \ dh \\
&= \frac{1}{|K|} \int_G \chi_{K}(t^{-1}gsh) \chi^1_{K}(h) \ dh
= f(t^{-1}gs). \end{align*}
Hence $(\chi^1_{sK} \cdot \chi_{tK})(g)\not=0$ only when
$g\in tKKs^{-1} = ts^{-1}KK$, that is, $st^{-1}g\in KK$.

For $g\in G$, we have that
\[ (a_n\cdot x)(g) = \frac1n \sum_{k=1}^\infty \sum_{l=1}^n
|K|^{-1} (\chi^1_{g_l K}\cdot\chi_{g_k K})(g)
= \frac1n \sum_{k=1}^\infty \sum_{l=1}^n f(g_k^{-1} g g_l). \]
Hence, for $g\not\in KK$, we see that
\[ \big|(a_n\cdot x)(g)\big| \leq \frac1n 
\big|\big\{ 1\leq l\leq n : g_l g_k^{-1} g \in KK
\text{ for some }k\geq1 \big\}\big| \leq \frac2n, \]
from the above.  If $b=(b_n)\in I_1$ then $\|b_n\|_\infty\rightarrow 0$,
and so $\| b_n|_{KK} \|_1\rightarrow0$, and from this it follows that
$ba \not= a$, exactly as in the proof of Theorem~\ref{dis_case}.
This contradiction completes the proof.
\end{proof}

Notice that, for abelian groups $G$, the above improves upon
Theorem~\ref{ab_com_case}, as Theorem~\ref{ab_com_case} only
tells us that $(L^1(G))_{\mc U}$ is not amenable when there
is a surjection $(L^1(G))_{\mc U} \rightarrow AP(G)'$; unless
$G$ is sufficiently ``small'', we cannot
necessarily take $\mc U$ to be an ultrafilter on a countable set.

We conclude with the following.

\begin{theorem}
Let $G$ be a locally compact group such that $L^1(G)$ is
ultra-amenable.  Then $G$ is finite; or $G$ satisfies the following:
\begin{enumerate}
\item $G$ is amenable;
\item $G$ is not compact nor a relative-[IN] group
   (so that $G$ is not abelian or discrete);
\item $AP(G)$ is finite-dimensional;
\item if $H$ is a closed normal subgroup of $G$ then either
   $G/H$ is finite, or $G/H$ satisfies the above properties.
\end{enumerate}
\end{theorem}

We currently do not know of any group which satisfies the above
conditions, so we strongly suspect that $L^1(G)$ is only ultra-amenable
when $G$ is finite.

We suspect that a careful argument using the ideas of
Proposition~\ref{quot_group} could reduce this problem to the
study of \emph{totally disconnected} groups (see \cite[Section~12.3]{Palmer2}),
as connected groups are fairly well understood (they are pro-Lie
groups, see \cite[Section~12.2]{Palmer2}).  Of course, totally
disconnected groups are not terribly well understood.  We strongly
suspect that the correct course of attack is to be improve
the proof of Theorem~\ref{in_case} so that it will hold for
all non-compact groups.

\section{Acknowledgments}

The results of Section~\ref{tensor_sec} are mostly from the author's
PhD thesis \cite{DawsThesis} completed at the University of
Leeds under the financial support of the EPSRC, and the guidance of
his PhD supervisors Garth Dales and Charles Read.  The author would
also like to thank George Willis, who proved to be an excellent source
for all questions group related, and Volker Runde, for useful
conversations, especially about C$^*$-algebras.  The author made great
use of the unpublished notes ``Ultra-methods'' by Mark Smith.  Finally,
the author wishes to thank the anonymous referee for careful proofreading.

\vspace{5ex}

\noindent\emph{Author's Address:}
\parbox[t]{1.4in}{St. John's College,\\
Oxford,\\
OX1 3JP.}

\bigskip\noindent\emph{Email:} \texttt{matt.daws@cantab.net}

\newpage

\appendix
\section{Errata}

\begin{abstract}
Some of the results of Section~5 of this paper are incorrect; in particular,
the characterisation of when an algebra is ultra-amenable, in terms of a diagonal
like construction, is not proved; and Theorem~5.7 is stated wrongly.
The rest of the paper is unaffected.
We shall show in this erratum that Theorem~5.7 can be corrected, and
that the other results of Section~5 are true if the algebra in
question has a certain approximation property.

\noindent\emph{2000 Mathematics subject classification:}
46B08; 46B28.
\end{abstract}

Some of the results of Section~5 of \cite{Daws1} are incorrect.  The claim
(ii)$\Rightarrow$(i) of Proposition~5.4 implicitly assumes that $\psi_0$ is bounded
below, but this is unproven.  Hence also the claim, in Corollary~5.5, that
if a Banach algebra $\mc A$ is contractible then it is ultra-amenable, is unproven.
Similarly (ii)$\Rightarrow$(i) of Theorem~5.6 requires $\psi_0$ to be bounded below.
The rest of the paper is unaffected.  We used some of these ideas in \cite[Section~4]{DR},
and so this is also incorrect; an erratum has been submitted.

Firstly, we deal with correcting Theorem~5.7.  We say that a C$^*$-algebra $\mc A$ is
\emph{subhomogeneous} if there exists $n\in\mathbb N$ such that every irreducible
representation of $\mc A$ has dimension at most $n$.  Subhomogeneous von Neumann
algebras have the special form claimed in Theorem~5.7, but this is not true for
C$^*$-algebras, see \cite[Section~IV.1.4]{blackadar} for examples.  This circle of ideas
was considered in \cite[Theorem~2.5]{llw} but we have been unable to follow
some of the proofs (in particular, the claim that (A4)$\Rightarrow$(R5)) so we provide
details here.

\begin{theorem}\label{cstar_ultra_amen}
Let $\mc A$ be a $C^*$-algebra.  Then the following are equivalent:
\begin{enumerate}
\item\label{cua:one} $\mc A$ is ultra-amenable;
\item\label{cua:two} $\mc A''$ is amenable;
\item\label{cua:three} $\ell^\infty(\mc A,I)$ is amenable for any index set $I$;
\item\label{cua:four} $\mc A$ is subhomogeneous.
\end{enumerate}
\end{theorem}
\begin{proof}
The original argument using the approximation property in \cite{Daws1}
is correct and shows (\ref{cua:one})$\Rightarrow$(\ref{cua:two}).
Similarly, as argued in \cite{Daws1} (see also \cite{llw}) if (\ref{cua:two})
holds, then $\mc A''$ has the form
\[ \mc A'' = \sum_{k=1}^n L^\infty(X_k) \otimes \mathbb M_{n_k}, \]
where for each $k$, $X_k$ is a measure space, and $n_k\in\mathbb N$.
Notice that if $\mc A''$ is of this form, then following \cite{llw},
it is elementary to see that so is $l^\infty(\mc A'',I)$ for any index set $I$.
This does imply that $\mc A$ is subhomogeneous
(see \cite[Proposition~IV.1.4.6]{blackadar}) but not that $\mc A$ has
the form originally claimed in \cite[Theorem~5.7]{Daws1}.

However, there is an algebraic characterisation of when $C^*$-algebras
are subhomogeneous, see \cite[Section~3.6]{dix}
or \cite[Section~IV.1.4.5]{blackadar}.  The algebra $\mathbb M_n$ is
of dimension $n^2$ and so for any $r>n^2$, we have
\begin{equation}
\sum_{\sigma\in S_r} \epsilon_\sigma x_{\sigma(1)} \cdots x_{\sigma(r)} = 0,
\label{polyr}
\end{equation}
for any $x_1,\cdots,x_r\in\mathbb M_n$ (this is readily seen by taking a basis).
Here $S_r$ is the symmetric group and $\epsilon:S_r\rightarrow\{\pm1\}$ the signature.
Let $r(n)$ be the smallest $r$ for which this holds for $\mathbb M_n$.
Then \cite[Lemma~3.6.2]{dix} shows that $r(n) \geq r(n-1)+2$
(see \cite[Section~IV.1.4.5]{blackadar} and references therein for better estimates).
As irreducible representations separate the points of a C$^*$-algebra $\mc A$, we conclude
that the following are equivalent:
\begin{enumerate}
\renewcommand{\labelenumi}{\roman{enumi}.}
\item any irreducible representation of $\mc A$ is of dimension at most $n$;
\item for any $x_1,\cdots,x_{r(n)}\in\mc A$, identity (\ref{polyr}) holds
for $r=r(n)$.
\end{enumerate}
Indeed, the only unclear issue is if $\pi:\mc A\rightarrow\mc B(H)$ is irreducible,
with $H$ infinite dimensional, why cannot (ii) hold?  However, then $\pi(\mc A)$ is
strongly dense in $\mc B(H)$, and $\mathbb M_{n+1}$ is a subalgebra
of $\mc B(H)$, which is enough to show that (ii) fails.

It is clear that the second condition passes to subalgebras, and with a little thought,
it is seen to pass to ultrapowers as well.  Thus, if (\ref{cua:two}) holds, then
$\ell^\infty(\mc A,I)$ is subhomogeneous and $\mc A$ is subhomogeneous, showing (\ref{cua:four}).
It is reasonably easy to show that $\ell^\infty(\mc A,I)$ is thus nuclear
(see \cite[Proposition~2.7.7]{BO}), or follow \cite[Theorem~2.5]{llw} for a
direct argument that $\ell^\infty(\mc A,I)$ is thus amenable.  As amenability passes
to quotients, (\ref{cua:three})$\Rightarrow$(\ref{cua:one}) is clear.
Finally, if (\ref{cua:four}) holds then any ultrapower of $\mc A$ is subhomogeneous
and hence amenable, showing (\ref{cua:one}).
\end{proof}

We erroneously claimed in \cite{DR} that (\ref{cua:one}) and
(\ref{cua:three}) are equivalent for any Banach algebra $\mc A$.  It would be
interesting to know if this is true.

We shall now improve \cite[Proposition~4.7]{Daws1}, and show that the map $\psi_0$ is
indeed bounded below for a wide class of Banach algebras $\mc A$.  We leave open
whether this holds for all $\mc A$ (which seems unlikely).  It seems possible
that similar, but stronger, conditions could characterise when $\mc A$ is ultra-amenable,
but we shall not pursue this here.

Let $E$ and $F$ be Banach spaces, and let $\mc U$ be an ultrafilter on an index set $I$.
As in \cite{Daws1}, we shall suppose that $\mc U$ is countably incomplete.
Recall from \cite[Section~4]{Daws1} the map $\psi_0:(E)_{\mc U}\proten(F)_{\mc U}
\rightarrow (E\proten F)_{\mc U}$, defined on elementary tensors by
\[ \psi_0( (x_i) \otimes (y_i) ) = (x_i\otimes y_i)
\qquad ( (x_i)\in (E)_{\mc U}, (y_i)\in (F)_{\mc U} ). \]
For the following, we recall that \cite[Theorem~9.1]{Hein1} characterises, in
terms of local properties, when an ultrapower has the (bounded) approximation
property.

\begin{theorem}
If $(E)_{\mc U}$ has the approximation property, then $\psi_0$ is an
injection for any $F$.
\end{theorem}
\begin{proof}
Let $\tau\in (E)_{\mc U}\proten(F)_{\mc U}$ have representation
$\tau = \sum_{n=1}^\infty x_n\otimes y_n$ with $\sum_n \|x_n\| \|y_n\|
< \infty$.
If $(E)_{\mc U}$ has the approximation property then, by \cite[Proposition~4.6]{Ryan1},
if $\tau\in (E)_{\mc U}\proten(F)_{\mc U}$ is non-zero, then there exist
$\mu\in (E)_{\mc U}'$ and $\lambda\in (F)_{\mc U}'$ with
\[ 0 \not= \ip{\mu\otimes\lambda}{\tau} =
\sum_{n=1}^\infty \ip{\mu}{x_n} \ip{\lambda}{y_n}. \]
As we only care about the value of $\mu$ on the countable set $\{x_n\}$,
by \cite[Corollary~7.5]{Hein1}, we may suppose that $\mu\in (E')_{\mc U}$,
and similarly, that $\lambda\in (F')_{\mc U}$, say $\mu=(\mu_i)$ and $\lambda=(\lambda_i)$.
Pick representatives $x_n = (x^{(i)}_n)$ and $y_n = (y^{(i)}_n)$, so that by absolute
convergence,
\[ \ip{\mu\otimes\lambda}{\tau} = \lim_{i\rightarrow\mc U} \sum_{n=1}^\infty
\ip{\mu_i}{x^{(i)}_n} \ip{\lambda_i}{y^{(i)}_n}
= \ip{(\mu_i\otimes\lambda_i)}{\psi_0(\tau)}. \]
Hence we must have that $\psi_0(\tau)\not=0$.
\end{proof}

Consequently, \cite[Corollary~5.5]{Daws1} correctly shows that if $\mc A$ is a
contractible Banach algebra with the approximation property, then $\mc A$ is
ultra-amenable.  However, a result of Selivanov, see \cite[Theorem~4.1.5]{RundeBook1},
shows that under these conditions, $\mc A$ is already the finite sum of full matrix algebras!

It is worth pointing out what can go wrong here (and hence the exact mistake in
the proof of \cite[Proposition~5.4]{Daws1}).  If $\mc A$ is contractible, then
we can find $\tau\in\mc A \proten \mc A$ with $a\cdot\tau=\tau\cdot a$ and
$\Delta(\tau)a=a$ for $a\in\mc A$.  We can then treat $\tau$ as a member of
$(\mc A)_{\mc U}\proten(\mc A)_{\mc U}$, and we have that $a\cdot\psi_0(\tau) =
\psi_0(\tau)\cdot a$ for $a\in(\mc A)_{\mc U}$.  As $\psi_0$ is an $(\mc A)_{\mc U}$-module
homomorphism, $\psi_0(a\cdot\tau-\tau\cdot a)=0$ for any $a\in(\mc A)_{\mc U}$.
However, if $\psi_0$ might fail to be injective, then this is not useful.

The following improves \cite[Proposition~4.7]{Daws1}, as a result of Grothendieck,
see \cite[Corollary~5.51]{Ryan1}, shows that a reflexive Banach space with the
approximation property automatically has the metric approximation property.

\begin{theorem}
If $(E)_{\mc U}$ has the bounded approximation property, then $\psi_0$ 
is bounded below.
\end{theorem}
\begin{proof}
Let $(E)_{\mc U}$ have the bounded approximation property with bound $M$, so by
(the obvious generalisation of) \cite[Theorem~4.14]{Ryan1}, the embedding
$(E)_{\mc U}\proten (F)_{\mc U}
\rightarrow \mc F((E)_{\mc U},(F)_{\mc U}')$ is bounded below by $M^{-1}$.
Here $\mc F((E)_{\mc U},(F)_{\mc U}')$ is the collection of finite-rank
operators $(E)_{\mc U}\rightarrow (F)_{\mc U}'$, given the operator norm.

Let $\tau\in (E)_{\mc U}\proten (F)_{\mc U}$ have representative
$\tau = \sum_{n=1}^\infty x_n \otimes y_n$.  For $\epsilon>0$, we can find
$T\in\mc F((E)_{\mc U},(F)_{\mc U}')$ with $\|T\|\leq M+\epsilon$ and
$|\ip{T}{\tau}| \geq \|\tau\|$.  Pick a representative
\[ T = \sum_{k=1}^N \mu_k \otimes \lambda_k, \]
for some $(\mu_k)\subseteq(E)_{\mc U}'$ and $(\lambda_k)\subseteq(F)_{\mc U}'$,
so that
\[ \ip{T}{\tau} = \sum_{n=1}^\infty \sum_{k=1}^N \ip{\mu_k}{x_n} \ip{\lambda_k}{y_n}. \]
Let $G$ be the closed span of $\{x_n\}$, so by \cite[Corollary~7.5]{Hein1}, we can find
a contraction $\phi:\lin\{\mu_k\}\rightarrow(E')_{\mc U}$ such that
\[ \ip{\phi(\mu_k)}{x} = \ip{\mu_k}{x} \qquad (1\leq k\leq N, x\in G). \]
It's not hard to see that then
\[ T_0 = \sum_{k=1}^N \phi(\mu_k)\otimes\lambda_k \]
satisfies $\|T_0\|\leq M+\epsilon$ and $\ip{T_0}{\tau} = \ip{T}{\tau}$.
In other words, we can assume that $\mu_k\in(E')_{\mc U}$ for each $k$;
analogously, we may also assume that $\lambda_k\in(F')_{\mc U}$ for each $k$.

So, pick representatives $\mu_k = (\mu^{(i)}_k)$ and $\lambda_k=(\lambda^{(i)}_k)$,
and for each $i$, let
\[ T_i = \sum_{k=1}^N \mu^{(i)}_k \otimes \lambda^{(i)}_k. \]
As $\mc U$ is countably incomplete, we can find a sequence $(\epsilon_i)$ of
strictly positive reals such that $\lim_{i\rightarrow\mc U} \epsilon_i = 0$.
For each $i$, pick $y_i\in E$ with $\|y_i\|\geq 1$ and $\|T_i(y_i)\|\geq \|T_i\|-\epsilon_i$.
Let $y=(y_i)$ so $\|y\|=1$ and $T(y) = (T_i(y_i))$ so that
\[ \lim_{i\rightarrow\mc U} \|T_i\| = \lim_{i\rightarrow\mc U} \|T_i\|-\epsilon_i
\leq \lim_{i\rightarrow\mc U} \|T_i(y_i)\| = \|T(y)\| \leq (M+\epsilon). \]
Finally, a calculation shows that
\[ \ip{T}{\tau} = \ip{(T_i)}{\psi_0(\tau)}, \]
where $(T_i) \in (\mc B(E,F^*))_{\mc U} \subseteq (E\proten F)_{\mc U}'$.
We conclude that $\|\psi_0(\tau)\| \geq \|\tau\| (M+\epsilon)^{-1}$, so that
$\psi_0$ is bounded below by $M^{-1}$.
\end{proof}

This shows that \cite[Theorem~5.6]{Daws1} does give a correct
characterisation of ultra-amenability for Banach algebras $\mc A$ whose
ultrapowers have the bounded approximation property.  This includes, for
example, algebras of the form $L^1(G)$ for a locally compact group $G$.

\medskip

\noindent\textbf{Acknowledgments:} I wish to thank: Seytek Tabaldyev who brought
the problems in \cite{DR} (and hence also in \cite{Daws1}) to my attention;
Volker Runde for pointing out the error in Theorem~5.7; and both Volker Runde
and the authors of \cite{llw} for suggesting to look at polynomial identities
as a way to fix Theorem~5.7.

\noindent
Matthew Daws,\\
School of Mathematics,\\
University of Leeds,\\
LEEDS LS2 9JT\\
United Kingdom\\
email: \texttt{matt.daws@cantab.net}

\end{document}